\documentclass[12pt,a4paper]{article}
\usepackage[dvips]{graphicx}
\usepackage{latexsym}
\usepackage{amssymb}
\usepackage{amsmath}
\usepackage{amscd}
\usepackage{mathrsfs}
\usepackage{bm}
\newtheorem{Th}{Theorem}[section]
\newtheorem{Co}[Th]{Corollary}
\newtheorem{Lem}[Th]{Lemma}

\newtheorem{Rem}[Th]{Remark}

\newtheorem{Pro}[Th]{Proposition}
\newcommand{\demo}{\par\noindent{\it Proof. \/}\ }
\newcommand{\enD}{\hfill $\Box$ \vspace{3truemm}\par}
\newcommand{\bx}{\mbox{\boldmath $x$}}
\newcommand{\bX}{\mbox{\boldmath $X$}}
\newcommand{\be}{\mbox{\boldmath $e$}}
\newcommand{\ba}{\mbox{\boldmath $a$}}
\newcommand{\bb}{\mbox{\boldmath $b$}}

\newcommand{\bv}{\mbox{\boldmath $v$}}
\newcommand{\by}{\mbox{\boldmath $y$}}
\newcommand{\bo}{\mbox{\boldmath $0$}}

\newcommand{\bxi}{\mbox{\boldmath $\xi$}}

\newcommand{\bt}{\mbox{\boldmath $t$}}

\newcommand{\bn}{\mbox{\boldmath $n$}}

\newcommand{\blambda}{\mbox{\boldmath $\lambda$}}

\newcommand{\R}{{\mathbb R}}
\newcommand{\lon}{\longrightarrow}

\newcommand{\sbxi}{\mbox{\scriptsize \boldmath$\xi$}}

\newcommand{\sblambda}{\mbox{\scriptsize \boldmath$\lambda$}}

\newcommand{\sbn}{\mbox{\scriptsize \boldmath$n$}}

\newcommand{\ou}{{\overline{u}}}
\begin{document}

\title{Caustics and Maxwell sets of world sheets \\ 
in anti-de Sitter space}

\author{Shyuichi IZUMIYA}

\date{\today}

\maketitle
\begin{abstract}
A world sheet in anti-de Sitter space is a timelike submanifold consisting of a one-parameter family of spacelike submanifolds.
We consider the family of lightlike hypersurfaces along spacelike submanifolds in the world sheet.
The locus of the singularities of lightlike hypersurfaces along spacelike submanifolds forms the caustic of the world sheet.
This notion is originally introduced by Bousso and Randall in theoretical physics. In this paper we give a mathematical framework
for the caustics of world sheets as an application of the theory of graph-like Legendrian unfoldings.
\end{abstract}

\section{Introduction}
In this paper we consider geometrical properties of caustics and Maxwell sets of world sheets in anti-de Sitter space
as an application of the theory of Legendrian unfoldings \cite{Izumiya93,Izumiya-Takahashi, Izumiya-Takahashi2,Izumiya-Takahashi3,Graph-like, GeomLag14} which is a special but an important case of the theory of wave front propagations \cite{Zak1}.
Anti-de Sitter space is one of the Lorentz space forms with rich geometric properties.
It is defined as a pseudo-sphere with a negative curvature in semi-Euclidean space with index 2 which admits the biggest symmetry in  Riemannian or Lorentz space forms.
Anti-de Sitter space plays important roles in theoretical physics such as the theory of general relativity, the 
string theory and the brane world scenario etc. It is one of the typical model of bulk spaces of the brane world scenario or the string theory (cf. \cite{Bousso,Bousso-Randall,KR01,M98,RS99,W98}). 
On the other hand, one of the important objects in the theoretical physics is the notion of lightlike hypersurfaces (light-sheets in physics) because they provide good models for  different types of horizons \cite{Ch, MTW}. In \cite{IzuLight} we considered lightlike hypersurfaces along spacelike submanifolds with general codimension in anti-de Sitter space.
lightlike hypersurfaces usually have singularities. 
We showed that lightlike hypersurfaces are wave fronts and applied the theory of Legendrian singularities \cite{Arnold1, Zak} to obtaining geometric properties of the singularities of lightlike hypersrufaces.
\par
 A world sheet (or a brane) in anti-de Sitter space is a timelike submanifold consisting of a one-parameter family of spacelike submanifolds.
Each spacelike submanifold is called a momentary space.
Since a momentary space is a spacelike submanifold, we have a lightlike hypersurface along each momentary space as a consequence of \cite{IzuLight}.
The set of singular values of a lightlike hypersurface is called the focal set along the momentary space. 
Since the world sheet is a one-parameter family of momentary spaces,
we naturally consider the family of lightlike hypersurfaces along momentary spaces in the world sheet.
The locus of the singularities (the focal sets) of lightlike hypersurfaces along momentary spaces is the caustic of the world sheet which was introduced by Bousso and Randall \cite{Bousso, Bousso-Randall} in order to define the notion of holographic domains. In this paper we construct a mathematical framework for the caustic of a world sheet and investigate the geometric properties of the singularities of the caustics of world sheets.
For the purpose, we apply the theory of graph-like Legendrian unfoldings \cite{Graph-like, GeomLag14}.
We also consider the notion of  Maxwell sets (crease sets) of world sheets which play an important role in the cosmology \cite{Penrose,Siino}.
In their paper \cite{Bousso, Bousso-Randall} the authors draw pictures on the simplest case (cf. \cite[Figures 2 and 3]{Bousso-Randall}).
However, this case the caustic coincides with the Maxwell set (i.e. a line).
In general, these sets are different, so that we consider both of them in this paper and emphasize that the Maxwell set of a world sheet is also an important subject.
\par
On the other hand, caustics appear in several area in physics (i.e. geometrical optics \cite{Nye}, the theory of underwater acoustics \cite{Brekhov} and the theory of gravitational lensings \cite{Peters}
, and so on) and mathematics (i.e. classical differential geometry \cite{Bru-Gib,Izu07,Porteous} and theory of differential equations \cite{Hormander,IzuHJ93}, and so on \cite{Arnold-pink}).
The notion of caustics originally belongs to geometrical optics. We can observe the caustic formed by the rays reflected at a mirror. One of the examples of caustics in the classical differential geometry is the evolute of a curve in the Euclidean 
plane which is given by the envelope of normal lines emanated from the curve.
The ray in the Euclidean plane is considered to be a line, so that the evolute is the caustic in the sense of geometrical optics.  Moreover, the singular points of the evolute correspond to the vertices of the original curve. The vertex is the point at where the curve has higher order contact with the osculating circle (i.e. the point where the curvature has an extremum).
Therefore, the evolute provides important geometrical information of the curve. We have the notion of evolutes for general hypersurfaces in the Euclidean space similar to
the plane curve case. In particular, there are detailed investigations on evolutes for surfaces in the Euclidean $3$-space \cite {Izu07,Porteous}.
Analogous to the Euclidean case, we can define the evolute of a hypersurface in Lorentz-Minkowski space \cite{Saloom,Farid}. Since a world sheet is a timelike submanifold, we may consider the evolute of
a timelike hypersurface in Lorentz-Minkowski space. However, the normal line is directed by a spacelike vector, so that the speed of the line exceeds the speed of the ray.
Although the evolute of a timelike hypersurface is a caustic in the theory of Lagrangian singularities, it is not a caustic in the sense of physics.
The situation in anti-de Sitter space is similar to that of Lorentz-Minkowski space. In a Lorentz manifold, the ray is directed by a lightlike vector, so that rays emanated from a spacelike submanifold
forms a lightlike hypersurface. Moreover, we have no notions of the time constant in the relativity theory. Hence everything that is moving depends on the time.
Therefore, we have to consider one parameter families of spacelike submanifolds (i.e. world sheets) in a Lorentz manifold, so that the notion of caustics by Bousso and Randall \cite{Bousso, Bousso-Randall} is essential. 
For further theoretical investigation, we construct a mathematical (geometric) framework for the caustics and the Maxwell sets of world sheets in this paper.
\par
We remark that the similar construction can be obtained for other Lorentz space forms (i.e. Lonrentz-Minkowski space and de Sitter space).
For a general Lorentz manifold, the situation is different from the case of Lorentz space forms.
In this case, we cannot construct explicit generating families for corresponding graph-like Legendrian unfoldings (cf. \S 6).
However, we can apply the theory of graph-like Legendrian unfoldings by using the classical method of characteristics for the (singular) eikonal equation corresponding to
the Lorentz metric. The detailed results will be appeared in elsewhere.

\section{Semi-Euclidean space with index 2}
\label{sec:1}
\par
In this section we prepare the basic notions on the semi-Euclidean
(n+2)-space with index 2. For detailed properties of the semi-Euclidean space,
see \cite{Oneil}.
For
any vectors $\bx=(x_{-1}, x_0, x_1,\cdots,x_{n}), \by
=(y_{-1}, y_0, y_1,\cdots,y_{n})\in \Bbb R^{n+2},$ the {\it
pseudo scalar product \/} of $\bx$ and $\by$ is defined to be
$\langle\bx,\by\rangle =-x_{-1} y_{-1}-x_0 y_0 +\sum_{i=1}^{n}x_i
y_i$. We call $(\Bbb R^{n+2}, \langle ,\rangle )$ a {\it
semi-Euclidean\/} (n+2)-{\it space with index 2\/} and write $\Bbb
R^{n+2}_2$ instead of $(\Bbb R^{n+2},\langle ,\rangle )$.
We say that a non-zero vector $\bx$ in $\Bbb R^{n+2}_2$ is {\it
spacelike\/}, {\it null\/} or {\it timelike\/} if
$\langle\bx,\bx\rangle>0,\langle\bx,\bx\rangle =0$ or
$\langle\bx,\bx\rangle <0$ respectively. The norm of the vector $\bx
\in \Bbb R^{n+2}_2$ is defined to be $\|\bx\|=\sqrt{|\langle\bx,
\bx\rangle|}$.
We define the {\it signature} of $\bx$ by
\[{\rm sign} (\bx)=\left\{
  \begin{array}{ccc}
  1\qquad\quad $\bx$\ \mbox{is\ spacelike}\\

  \mbox{}0\qquad\quad \bx\ \mbox{is\ null} \hspace*{\fill}\\

  -1\qquad \bx\ \mbox{is\ timelike}
  \end{array}\right.
  \]
 For a non-zero vector $\bn\in \Bbb R^{n+2}_2$ and a real number $c$, we define a
 {\it hyperplane with pseudo-normal \/}$\bn$ by
$$
HP(\bn,c)=\{\bx \in \Bbb R^{n+2}_2 |\langle\bx,\bn\rangle=c\}.
$$
We call $HP(\bn,c)$ a {\it Lorentz hyperplane\/}, a {\it
semi-Euclidean hyperplane with index 2\/} or a {\it null
hyperplane\/} if $\bn$ is {\it timelike, spacelike or null
\/}respectively.

We now define the {\it Anti de Sitter $n+1$-space \/} (briefly, the {\it AdS
$n+1$-space\/}) by
$$
AdS^{n+1}=\{\bx\in \Bbb R_2^{n+2}\ |\ \langle\bx,\bx\rangle =-1\}=H^{n+1}_1,
$$
the {\it unit pseudo $n+1$-sphere with index 2\/} by
$$
S^{n+1}_2=\{\bx\in \R_2^{n+2}\ |\ \langle\bx,\bx\rangle =1\},
$$
and the {\it {\rm (}closed{\rm )} nullcone\/} with vertex $\bm{\lambda}\in \R^{n+2}_2$ by\\
$$\Lambda_{\bm{\lambda}}^{n+1}=\{\bx\in \R_{2}^{n+2}|\langle\bx-\bm{\lambda}, \bx-\bm{\lambda}\rangle=0\}.
$$
In particular we write $\Lambda ^*=\Lambda ^{n+1}_0\setminus \{\bo\}$ and also call it
the {\it {\rm (}open{\rm )} nullcone\/}. Our main subject in this paper is $AdS^{n+1}$.
Since the causality of $AdS^{n+1}$ is violated, it is usually considered the universal covering
space $\widetilde{AdS}^{n+1}$ of $AdS^{n+1}$ in physics which is called the {\it universal Anti de Sitter space\/}.
We remark that the local structure of these spaces are the same.
Since $AdS^{n+1}$ is a Lorentz space form, there exists a lightcone on each tangent space.
Such a lightcone is explicitly expressed as follows:
For any $\bm{\lambda} \in AdS^{n+1},$ we have a hyperplane $HP(\bm{\lambda},-1).$
This hyperplane is the tangent hyperplane of $AdS^{n+1}$ at $\bm{\lambda}.$
We can show that 
\[
HP(\bm{\lambda},-1)\cap AdS^{n+1}=\Lambda ^{n+1}_{\bm{\lambda}}\cap AdS^{n+1}.
\]
Therefore, $HP(\bm{\lambda},-1)\cap AdS^{n+1}$ is the lightcone in the tangent hyperplane $HP(\bm{\lambda},-1)$
of $AdS^{n+1}$ at $\bm{\lambda}.$
We write it by $LC^{AdS}(\bm{\lambda})$ and call an {\it anti-de Sitter lightcone} (briefly, an {\it AdS-lightcone})
at $\bm{\lambda}\in AdS^{n+1}$.
\par
For any $\bx_1,\cdots, \bx_{n+1} \in \R^{n+2}_2$, we define a
vector $\bx_1\wedge \cdots\wedge \bx_n$ by
$$
\bx_1\wedge\cdots\wedge \bx_{n+1}= \vmatrix
-\be_{-1}&-\be_0&\be_1&\cdots&\be_{n}\vspace{2mm}\\
x^1_{-1}&x^1_0&x^1_1&\cdots&x^1_{n}\\
\vdots&\vdots&\vdots&\vdots&\vdots\\
x^{n+1}_{-1}&x^{n+1}_{0}&x^{n+1}_1&\cdots&x^{n+1}_{n}
\endvmatrix,
$$
where $\{\be_{-1}, \be_0, \be_1,\cdots,\be_{n}\}$ is the canonical
basis of $\Bbb R^{n+2}_2$ and $\bx_i=(x^i_{-1}, x^i_0,
x^i_1,\cdots,x^i_{n})$. We can easily check that
$$
\langle\bx,\ \bx_1\wedge\cdots\wedge
\bx_{n+1}\rangle=\textrm{det}(\bx, \bx_1,\cdots,\bx_{n+1}),
$$
so that $\bx_1\wedge\cdots\wedge \bx_{n}$ is pseudo-orthogonal to
any $\bx_i$\ (for\ $i=1,\cdots,n$).

\section{World sheets in in anti-de Sitter space}
\par
In this section we introduce the basic geometrical framework for the 
study of world sheets in anti-de Sitter $n+1$-space.
Consider the orientation of $\R^{n+2}_2$
provided by the condition that $\textrm{det}(\be_{-1}, \be_0,\be _1,\cdots,\be_{n})>0.$
This orientation induces the orientation of $x_{-1}x_0$-plane,
so that it gives a time
orientation on $AdS^{n+1}$. 
If we consider the universal Anti de Sitter space $\widetilde{AdS}^{n+1},$ we can determine the
future direction.
The world sheet is defined to be a timelike submanifold foliated by
a codimension one spacelike submanifolds.
Here, we only consider the local situation, so that we considered a one-parameter family of spacelike submanifolds.
Let $AdS^{n+1}$ be the oriented and time-oriented
anti-de Sitter space. 
Let $\bX :U\times I\lon AdS^{n+1}$ be a timelike embedding of codimension $k-1,$
where $U\subset \R^s$ ($s+k=n+2$) is an open subset and $I$ an open interval.
 We write 
$W=\bX(U\times I)
$ and identify $W$ and $U\times I$ through the embedding $\bX.$
Here, the embedding $\bX$ is said to be {\it timelike} if the tangent space $T_p  W$
of $W$ at $p=\bX(u,t)$ is a timelike subspace (i.e., Lorentz subspace of $T_pAdS^{n+1}$) for any point $p\in W$.
We write $\mathcal{S}_t=\bX(U\times\{t\})$ for each $t\in I.$
We call $\mathcal{S}=\{\mathcal{S}_t\ |t\in I\}$ a {\it spacelike foliation} 
on $W$ if $\mathcal{S}_t$ is a spacelike submanifold for any $t\in I.$
Here, we say that $\mathcal{S}_t$ is {\it spacelike} if the tangent space $T_p\mathcal{S}_t$ 
consists only spacelike vectors (i.e., spacelike subspace) for
any point $p\in \mathcal{S}_t.$
We call $\mathcal{S}_t$ a {\it momentary space} of $\mathcal{S}=\{\mathcal{S}_t\ |t\in I\}$.
For any $p=\bX(u,t)\in W\subset AdS^{n+1},$
we have
\[
T_pW=\langle \bX_t(u,t),\bX_{u_1}(u,t),\dots ,\bX_{u_s}(u,t)\rangle _\R,
\]
where $\bX_t=\partial \bX/\partial t,\bX_{u_j}=\partial \bX/\partial u_j.$
We say that $(W,\mathcal{S})$ (or, $\bX$ itself) is a {\it world sheet}
if $W$ is time-orientable.
Since $W$ is time-orientable, there exists a timelike vector field $\bv(u,t)$ on $W$ \cite[Lemma 32]{Oneil}.
Moreover, we can choose that $\bv$ is adapted with respected to
the time-orientation of $AdS^{n+1}.$
Here, we say that a timelike vector field $\bv(u,t)$ on $W$ is {\it adapted}
if 
$
\textrm{det}(\bX(u,t),\bv(u,t),\be _1,\dots ,\be _{n})>0.
$
Let $N_p(W)$ be the pseudo-normal space of $W$ at $p=\bX(u,t)$ in $\R^{n+2}_2.$
Since $T_pW$ is a timelike subspace of $T_p\R^{n+2}_2,$
$N_p(W)$ is a $k$-dimensional Lorentz subspace of $T_p\R^{n+2}_2$.
(cf.,\cite{Oneil}).
On the pseudo-normal space $N_p(W),$ we have a $(k-1)$-dimensional spacelike subspace:
\[
N^{AdS}_p(W)= \{\bxi\in N_p(W)\ |\ \langle \bxi,\bX(u,t)\rangle =0\ \},
\]
so that we have a $(k-2)$-unit sphere
\[
N^{AdS}_1(W)_p=\{\bxi\in N^{AdS}_p(W)\ |\ \langle \bxi,\bxi\rangle =1\ \}.
\]
Therefore, we have a unit spherical normal bundle over $W$:
\[
N^{AdS}_1(W)=\bigcup _{p\in W} N^{AdS}_1(W)_p.
\]
\par
On the other hand,
we write $N_p(\mathcal{S}_t)$ as the pseudo-normal space of $\mathcal{S}_t$ at $p=\bX(u,t)$
in $\R^{n+2}_2.$
Then $N_p(\mathcal{S}_t)$ is a $k+1$-dimensional semi-Euclidean subspace with index $2$ of $T_p\R^{n+2}_2$ \cite{Oneil}.
On the pseudo-normal space $N_p(\mathcal{S}_t),$ we have two kinds of pseudo spheres:
\begin{eqnarray*}
N_p(\mathcal{S}_t;-1)& = & \{\bv\in N_p(\mathcal{S}_t)\ |\ \langle \bv,\bv\rangle =-1\ \} \\
N_p(\mathcal{S}_t;1)&= & \{\bv\in N_p(\mathcal{S}_t)\ |\ \langle \bv,\bv\rangle =1\ \}.
\end{eqnarray*}
We remark that $N_p(\mathcal{S}_t;-1)$ is the $k$-dimensional anti-de Sitter space and $N_p(\mathcal{S}_t;1)$ is the $k$-dimensional pseudo-sphere with index $2.$
Therefore,
we have two unit spherical normal bundles $N(\mathcal{S}_t;-1)$ and $N(\mathcal{S}_t;1)$ over $\mathcal{S}_t$.
By definition, $\bX(u,t)$ is one of the timelike unit normal vectors of $\mathcal{S}_t$ at $p=\bX(u,t),$ so that
$\bX(u,t) \in N_p(\mathcal{S}_t).$ 
Since $\mathcal{S}_t=\bX(U\times \{t\})$ is a codimension one spacelike submanifold in $W,$ there exists a
unique timelike adopted unit normal vector field $\bn^T(u,t)$ of $\mathcal{S}_t$ such that 
$\bn^T(u,t)$ is tangent to $W$ at any point $p=\bX(u,t).$
It means that $\bn^T(u,t)\in N_p(\mathcal{S}_t)\cap T_pW$ with $\langle \bn^T(u,t),\bn^T(u,t)\rangle =-1$ and
$\det (\bX(u,t),\bn^T(u,t),\be _1,\dots ,\be_n ) >0.$
We define a $(k-2)$-dimensional spacelike unit sphere in $N_p(\mathcal{S}_t)$ by
\[
N^{AdS}_1(\mathcal{S}_t)_p[\bn ^T]=\{\bxi \in N_p(\mathcal{S}_t;1)\ |\ \langle \bxi, \bn ^T(u,t)\rangle =\langle \bxi,\bX(u,t)\rangle=0,p=\bX(u,t)\ \}.
\]
Then we have a {\it spacelike unit $(k-2)$-spherical bundle $N_1(\mathcal{S}_t)[\bn ^T]$ over $\mathcal{S}_t$ with respect to $\bn ^T$}.
Since we have
$T_{(p,\xi)}N^{AdS}_1(\mathcal{S}_t)[\bn^T]=T_p\mathcal{S}_t\times T_\xi N^{AdS}_1(\mathcal{S}_t)_p[\bn ^T],$
we have the canonical Riemannian metric on $N^{AdS}_1(\mathcal{S}_t)[\bn^T]$
which we write $(G_{ij}((u,t),\bxi))_{1\leqslant i,j\leqslant n-1}.$
Since $\bn^T$ is uniquely determined, we can write $N_1^{AdS}[\mathcal{S}_t]=N_1^{AdS}(\mathcal{S}_t)[\bn^T].$
Moreover, we remark that $N_1^{AdS}(W)|\mathcal{S}_t=N_1^{AdS}[\mathcal{S}_t]$ for any $t\in I.$
\par
We now define a map $\mathbb{NG}:N^{AdS}_1(W)\lon \Lambda ^*$ by $\mathbb{NG}(\bX(u,t),\bxi)=\bn^T(u,t)+\bxi$.
We call $\mathbb{NG}$ an {\it $AdS$-world nullcone Gauss image} of
$W=\bX(U\times I)$.
A {\it momentary nullcone Gauss image} of $N_1^{AdS}[\mathcal{S}_t]$ is defined 
to be the restriction of the $AdS$-world nullcone Gauss image
\[
\mathbb{NG}(\mathcal{S}_t)=\mathbb{NG}|N_1^{AdS}[\mathcal{S}_t]:N_1^{AdS}[\mathcal{S}_t]\lon \Lambda ^*.
\]
This map leads us to the notions of curvatures.
Let $T_{(p,\xi)}N_1[\mathcal{S}_t]$ be the tangent space of $N_1[\mathcal{S}_t]$ at $(p,\bxi).$
Under the canonical identification $(\mathbb{NG}(\mathcal{S}_t)^*T\R^{n+2}_2)_{(p,\sbxi)}
=T_{(\sbn^T(p)+\sbxi)}\R^{n+1}_1\equiv T_p\R^{n+2}_2,$
we have
\[
T_{(p,\sbxi)}N_1[\mathcal{S}_t]=T_p\mathcal{S}_t\oplus T_\xi S^{k-2}\subset T_pM\oplus N_p(\mathcal{S}_t)=T_p\R^{n+2}_2,
\]  
where $T_\xi S^{k-2}\subset T_\xi N_p(\mathcal{S}_t)\equiv N_p(\mathcal{S}_t)$ and $p=\bX(u,t).$
Let 
\[
\Pi ^t :\mathbb{NG}(\mathcal{S}_t)^*T\R^{n+2}_2=TN_1[\mathcal{S}_t]\oplus \R^{k+1}
\lon TN_1[\mathcal{S}_t]
\]
be the canonical projection.
Then
we have
a linear transformation
\[
S_N (\mathcal{S}_t)_{(p,\sbxi)}=-\Pi^t_{\mathbb{NG}(\mathcal{S}_t)(p,\xi)}\circ d_{(p,\xi)}\mathbb{NG}(\mathcal{S}_t)
: T_{(p,\xi)}N^{AdS}_1[\mathcal{S}_t]\lon T_{(p,\xi)}N^{AdS}_1[\mathcal{S}_t],
\]
which is called a {\it momentary nullcone shape operator} of $N^{AdS}_1[\mathcal{S}_t]$ at $(p,\bxi).$ 
\par
On the other hand, we choose a pseudo-normal section 
$\bn ^S(u,t)\in N^{AdS}_1(W)$ 
at least locally. Then we have 
$\langle \bn ^S,\bn ^S\rangle =1$ and $\langle \bX_t,\bn^S\rangle = \langle \bX_{u_i},\bn^S\rangle =\langle \bn ^T,\bn ^S\rangle =0,$
so that the vector
$\bn^T (u,t)+ \bn^S(u,t)$ is lightlike. 
We define a mapping
\[
\mathbb{NG}(\mathcal{S}_{t_0};\bn^S):U\lon \Lambda ^*
\]
by $\mathbb{NG}(\mathcal{S}_{t_0};\bn^S)(u)=\bn^T(u,t_0)+\bn^S(u,t_0),$
which is called a {\it momentary nullcone Gauss images of 
$\mathcal{S}_{t_0}=\bX(U\times \{t_0\})$ with respect to
$\bn^S.$}
Under the identification of $\mathcal{S}_{t_0}$ and $U\times\{t_0\}$ through $\bX,$ we have the
linear mapping provided by the derivative of the momentary nullcone Gauss image $\mathbb{NG}(\mathcal{S}_{t_0};\bn^S)$ at each point 
$p=\bX(u,t_0)$,
\[
d_p\mathbb{NG}(\mathcal{S}_{t_0};\bn^S):T_p\mathcal{S}_{t_0}\lon T_p\R^{n+1}_1= T_p\mathcal{S}_{t_0}\oplus N_p(\mathcal{S}_{t_0}).
\]
Consider the orthogonal projection $\pi ^t:T_p\mathcal{S}_{t_0}\oplus
N_p(\mathcal{S}_{t_0})\rightarrow T_p\mathcal{S}_{t_0}.$ We define
\[
S_p(\mathcal{S}_{t_0};\bn^S)=-\pi^t\circ d_p\mathbb{NG}(\mathcal{S}_{t_0};\bn^S):T_p\mathcal{S}_{t_0}\lon T_p\mathcal{S}_{t_0}.
\]
We call the
linear transformation $S_{p}(\mathcal{S}_{t_0};\bn^S)$ a {\it momentary $\bn^S$-shape
operator} of $\mathcal{S}_{t_0}=\bX (U\times \{t_0\})$ at $p=\bX (u,t_0).$ 
Let $\{\kappa
_{i}(\mathcal{S}_{t_0};\bn^S)(p)\}_ {i=1}^s$ be the eigenvalues of $S_{p}(\mathcal{S}_{t_0};\bn^S)$, which are called  {\it 
momentary nullcone
principal curvatures of $\mathcal{S}_{t_0}$ with respect to $\bn^S$\/} at $p=\bX(u,t_0)$.
Then a {\it momentary nullcone Gauss-Kronecker curvature of $\mathcal{S}_{t_0}$ with respect to
$\bn^S$\/} at $p=\bX (u,t_0)$ is defined to be
\[
K_N(\mathcal{S}_{t_0};\bn^S)(p)={\rm det} S_{p}(\mathcal{S}_{t_0};\bn^S).
\]
We say that a point $p=\bX (u,t_0)$ is a {\it momentary $\bn^S$-nullcone umbilical
point} of $\mathcal{S}_{t_0}$ if 
\[
S_{p}(\mathcal{S}_{t_0};\bn^S)=\kappa (\mathcal{S}_{t_0};\bn^S)(p) 1_{T_{p}\mathcal{S}_{t_0}}.
\]
We say that $W=\bX (U\times I)$ is {\it totally
$\bn^S$-nullcone umbilical} if any point $p=\bX(u,t)\in W$ is
momentary $\bn^S$-nullcone umbilical.
Moreover, $W=\bX(U\times I)$ is said to be {\it totally nullcone umbilical} if
it is totally $\bn^S$-nullcone umbilical for any $\bn^S.$
We deduce now the nullcone Weingarten formula. Since $\bX _{u_i}$
$(i=1,\dots s)$ are spacelike vectors, we have a Riemannian metric
(the {\it first fundamental form \/}) on $\mathcal{S}_{t_0}=\bX (U\times\{t_0\})$
defined by $ds^2 =\sum _{i=1}^{s} g_{ij}du_idu_j$,  where
$g_{ij}(u,t_0) =\langle \bX _{u_i}(u,t_0 ),\bX _{u_j}(u,t_0)\rangle$ for any
$u\in U.$ We also have a {\it nullcone second fundamental invariant of $\mathcal{S}_{t_0}$
with respect to the normal vector field $\bn ^S $\/} defined
by $h _{ij}(\mathcal{S}_{t_0};\bn^S )(u,t_0)=\langle -(\bn^T +\bn^S)
_{u_i}(u,t_0),\bX_{u_j}(u,t_0)\rangle$ for any $u\in U.$
By the similar arguments to those in the proof of \cite[Proposition 3.2]{IzuSM}, we have 
the following proposition.
\begin{Pro}
Let $\{\bX, \bn^T,\bn^S_1,\dots ,\bn^S_{k-1}\}$ be a a pseudo-orthonormal frame of $N(\mathcal{S}_{t_0})$ with $\bn^S_{k-1}=\bn^S.$ Then we have the following momentary nullcone Weingarten formulae {\rm :}
\vskip1.5pt
\par\noindent
{\rm (a)} $\mathbb{NG}(\mathcal{S}_{t_0};\bn^S)_{u_i}=\langle \bn ^T_{u_i},\bn ^S\rangle(\bn^T+\bn^S)+\sum _{\ell =1}^{k-2}\langle (\bn^T+\bn^S)_{u_i},\bn^S_\ell \rangle\bn^S_\ell -\sum_{j=1}^{s}
h_i^j(\mathcal{S}_{t_0};\bn^S )\bX _{u_j}$
\par\noindent
{\rm (b)} $
\pi ^t\circ \mathbb{NG}(\mathcal{S}_{t_0};\bn^S)_{u_i}=-\sum_{j=1}^{s}
h_i^j(\mathcal{S}_{t_0};\bn^S )\bX _{u_j}.
$
\smallskip
\par\noindent
Here $\displaystyle{\left(h_i^j(\mathcal{S}_{t_0};\bn^S )\right)=\left(h_{ik}(\mathcal{S}_{t_0};\bn^S)\right)\left(g^{kj}\right)}$
and $\displaystyle{\left( g^{kj}\right)=\left(g_{kj}\right)^{-1}}.$
\end{Pro}
\par
Since $\mathbb{NG}(\mathcal{S}_{t_0};\bn^S)_{u_i}=d\mathbb{NG}(\mathcal{S}_{t_0};\bn^S)(\bX_{u_i})$,
we have 
$$
S_p(\mathcal{S}_{t_0};\bn^S )(\bX_{u_i}(u,t_0))=-\pi^t\circ \mathbb{NG}(\mathcal{S}_{t_0};\bn^S)_{u_i}(u,t_0),
$$
so that the representation matrix of $S_p(\mathcal{S}_{t_0};\bn^S )$ with respect to the basis 
\[
\{\bX_{u_1}(u,t_0),\bX_{u_2}(u,t_0),\dots ,\bX_{u_s}(u,t_0)\}
\]
 of $T_p\mathcal{S}_{t_0}$ is $(h^i_j(\mathcal{S}_{t_0};\bn^S)(u,t_0)).$
Therefore, we have an explicit
expression of the momentary nullcone Gauss-Kronecker curvature of $\mathcal{S}_{t_0}$ with respect to
$\bn^S$ by
$$
K_N (\mathcal{S}_{t_0};\bn^S )(u,t_0)=\frac{\displaystyle{{\rm det}\left(h_{ij}(\mathcal{S}_{t_0};\bn^S )(u,t_0)\right)}}
{\displaystyle{{\rm det}\left(g_{\alpha \beta}(u,t_0)\right)}}.
$$
Since $\langle -(\bn^T +\bn^S )(u,t),\bX _{u_j}(u,t)\rangle =0$  we have
\[
h_{ij}(\mathcal{S}_{t_0};\bn^S)(u,t)=\langle \bn^T (u,t)+\bn^S (u,t),\bX
_{u_iu_j}(u,t)\rangle.
\]
 Therefore the momentary nullcone second fundamental
invariant of $\mathcal{S}_{t_0}$ at a point $p_0=\bX (u_0,t_0)$ depends only on the values 
$\bn^T (u_0)+\bn^S (u_0)$ and $\bX _{u_iu_j}(u_0)$, respectively.
Therefore, we write
\[
h_{ij}(\mathcal{S}_{t_0};\bn^S)(u_0,t_0)=h_{ij}(\mathcal{S}_{t_0})(p_0,\bxi_0),
\]
where $p_0=\bX(u_0,t_0)$ and $\bxi_0=\bn^S(u_0,t_0)\in N^{AdS}_1(W)_{p_0}.$
Thus, the momentary $\bn^S$-shape operator and the momentary nullcone curvatures also depend only on
$\bn^T (u_0,t_0)+\bn^S (u_0,t_0)$, $\bX_{u_i}(u_0,t_0)$  and $\bX
_{u_iu_j}(u_0,t_0)$, independent of the derivation of the vector fields
$\bn^T$ and 
$\bn^S .$ 
We may write $S_{p_0}(\mathcal{S}_{t_0};\bxi_0)=S_{p_0}(\mathcal{S}_{t_0};\bn^S),$ $\kappa _i(\mathcal{S}_{t_0},\bxi_0)(p_0)= \kappa _i(\mathcal{S}_{t_0};\bn^S)(p_0)$ $(i=1,\dots ,s)$
and $K_N(\mathcal{S}_{t_0},\bxi_0)(p_0)=K_N (\mathcal{S}_{t_0};\bn^S)(p_0)$ at $p_0=\bX (u_0,t_0)$
with respect to $\bxi_0=\bn^S(u_0,t_0).$ 
We also say that a point $p_0=\bX (u_0,t_0)$ is 
{\it momentary $\bxi_0$-nullcone umbilical\/} if $S_{p_0}(\mathcal{S}_{t_0};\bxi_0)=\kappa _i(\mathcal{S}_{t_0})(p_0,\bxi_0)1_{T_{p_0}\mathcal{S}_{t_0}}$. 
The momentary space $\mathcal{S}_{t_0}$ is said to be {\it totally momentary nullcone umbilical} if any point 
$p=\bX(u,t_0)$ is momentary $\bxi$-nullcone umbilical for any $\bxi\in N^{AdS}_1(\mathcal{S}_{t_0})_p[\bn^T]$.
Moreover, we say that a point $p_0=\bX (u_0,t_0)$ is a {\it momentary $\bxi_0$-nullcone parabolic point \/} of $W$ if
$K_N (\mathcal{S}_{t_0};\bxi_0)(p_0)=0.$ 
Let $\kappa _N(\mathcal{S}_{t})_i(p,\bxi)$ be the eigenvalues of the momentary nullcone shape operator $S_N(\mathcal{S}_{t}) _{(p,\sbxi)}$, $(i=1,\dots ,n-1)$. 
We write $\kappa _N(\mathcal{S}_t)_i(p,\bxi)$, $(i=1,\dots ,s)$ as the eigenvalues belonging to
the eigenvectors on $T_p\mathcal{S}_t$
and $\kappa _N(\mathcal{S}_t)_i(p,\bxi)$, $(i=s+1,\dots n)$ as the eigenvalues belonging to the eigenvectors on 
the tangent space of the fiber  
of $N_1[\mathcal{S}_t].$  
\begin{Pro}
For $p_0=\bX(u_0,t_0)$ and $\bxi_0\in N^{AdS}_1[\mathcal{S}_{t_0}]_{p_0},$ we have
$$\kappa _N(\mathcal{S}_{t_0})_i(p_0,\bxi_0)=\kappa _i(\mathcal{S}_{t_0},\bxi_0)(p_0),\ (i=1,\dots s),\ \kappa _N(\mathcal{S}_{t_0})_i(p_0,\bxi_0)=-1,\ (i=s+1,\dots n).$$
\end{Pro}
We call $\kappa _N(\mathcal{S}_{t})_i(p,\bxi)=\kappa _i(\mathcal{S}_{t},\bxi)(p)$, $(i=1,\dots ,s)$
the {\it nullcone principal curvatures} of $\mathcal{S}_{t}$ with respect to $\bxi$ at $p=\bX(u,t)\in W.$
\demo
Since $\{\bX, \bn^T,\bn^S_1,\dots ,\bn^S_{k-1}\}$ is a pseudo-orthonormal frame of $N(\mathcal{S}_t)$
and 
$$
\bxi_0=\bn^S_{k-1}(\ou_0,t_0)\in S^{k-2}=N_1[\mathcal{S}_{t_0}]_p,
$$ 
we have $
\langle \bn^T(\ou_0,t_0),\bxi_0\rangle =\langle \bn^S_i(\ou_0,t_0),\bxi_0\rangle =0$
for $i=1,\dots ,k-2.$
Therefore, we have 
$$T_{\bm{\xi}_0}S^{k-2}=\langle \bn^S_1(\ou_0,t_0),\dots ,\bn^S_{k-2}(\ou_0,t_0)\rangle .$$
By this orthonormal basis of $T_{\sbxi_0}S^{k-2},$
the canonical Riemannian metric $G_{ij}(p_0,\bxi_0)$ is represented by
\[
(G_{ij}(p_0,\bxi_0))=\left(
\begin{array}{cc}
g_{ij}(p_0)  & 0 \\
0 & I_{k-2}
\end{array}
\right) ,
\]
where $g_{ij}(p_0)=\langle \bX_{u_i}(\ou_0,t_0), \bX_{u_j}(\ou_0,t_0)\rangle $.
\par
On the other hand, by Proposition 3.1, we have
\[
-\sum_{j=1}^s h^j_i(\mathcal{S}_{t_0},\bn^S)\bX_{u_j}=\mathbb{NG}(\mathcal{S}_{t_0},\bn^S)_{u_i}=
d_{p_0}\mathbb{NG}(\mathcal{S}_{t_0};\bn^S)\left(\frac{\partial}{\partial u_i}\right),
\]
so that we have
\[
S_{p_0}(\mathcal{S}_{t_0};\bxi_0)\left(\frac{\partial}{\partial u_i}\right)=\sum_{j=1}^s h^j_i(\mathcal{S}_{t_0};\bn^S)\bX_{u_j}.
\]
Therefore, the representation matrix of $S_{p_0}(\mathcal{S}_{t_0};\bxi_0)$ with respect to the basis
$$
\{\bX_{u_1}(\ou_0,t_0),\dots ,\bX_{u_s}(\ou_0,t_0),\bn^S_1(\ou_0,t_0),\dots ,\bn^S_{k-2}(\ou_0,t_0)\}
$$ of $T_{(p_0,\bm{\xi}_0)}N_1[\mathcal{S}_{t_0}]$
is of the form
\[
\left(
\begin{array}{cc}
h^j_i(\mathcal{S}_{t_0},\bn^S)(u_0,t_0)  & * \\
0 & -I_{k-2}
\end{array}
\right).
\]
Thus, the eigenvalues of this matrix are $\lambda _i=\kappa _i(\mathcal{S}_{t_0},\bxi_0)(p_0)$, $(i=1,\dots ,s)$ and
$\lambda _i=-1,$ $(i=s+1,\dots ,n-1)$.
This completes the proof.

\enD

\section{Lightlike hypersurfaces along momentary spaces}
We define a hypersurface
$
\mathbb{LH}_{\mathcal{S}_t}:N^{AdS}_1[\mathcal{S}_t]\times \R\lon AdS^{n+1}
$
by
$$
\mathbb{LH}_{\mathcal{S}_t}(((u,t),\bxi),\mu)=\bX(u,t)+\mu (\bn^T(u,t)+\bxi)=\bX(u,t)+\mu\mathbb{NG}(\mathcal{S}_t)((u,t),\bxi),
$$
where $p=\bX (u,t),$ which is called a {\it momentary lightlike hypersruface\/} in anti-de Sitter space along $\mathcal{S}_t$.
We remark that $\mathbb{LH}_{\mathcal{S}_t}(N^{AdS}_1[\mathcal{S}_t]\times\R)$ is a lightlike hypersurface.
Here a hypersurface is {\it lightlike} if the tangent space of the hypersurface at any regular point is a lightlike hyperplane.
 \par
We define a family of functions $H: U\times I\times AdS^{n+1}\lon \R$ on a world sheet  $W=\bX (U\times I)$ 
 by
$
 H((u,t),\blambda )=\langle \bX (u,t) ,\blambda\rangle +1.
$
We call
$H$ the {\it anti-de Sitter height function\/} (briefly, AdS-height function) on the world sheet
$W=\bX(U\times I).$
For any fixed $(t_0,\blambda _0)\in I\times\R_2^{n+2},$ we write $h_{(t_0,\sblambda _0)}(u)=H((u,t_0),\blambda _0).$
\begin{Pro}
Let $W$ be a world sheet 
and
$H: U\times I\times(AdS^{n+1}\setminus W)\to\R$
the AdS-height function on $W.$
Suppose that $p_0=\bX(u_0,t_0)\not=\blambda _0.$ Then we have the following$:$
\par\noindent
{\rm (1)}
$h_{(t_0,\sblambda _0)}(u_0)=\partial h_{(t_0,\sblambda _0)}/\partial u_i(u_0)=0$, $(i=1,\dots ,s)$
if and only if
there exist $\bxi_0 \in N^{AdS}_1[\mathcal{S}_{t_0}]_{p_0}$ and $\mu_0\in
\R\setminus \{0\}$ such that 
$
\blambda _0 =\mathbb{LH}_{\mathcal{S}_{t_0}}(((u_0,t_0),\bxi_0),\mu_0).
$
\par\noindent
{\rm (2)}
$h_{(t_0,\sblambda _0)}(u_0)=\partial h_{(t_0,\sblambda _0)}/\partial u_i(u_0)=
{\rm det}{\mathcal H}(h_{(t_0,\sblambda _0)})(u_0)=0$ $(i=1,\dots ,s)$
if and only if
there exist $\bxi_0 \in N_1[\mathcal{S}_{t_0}]_{p_0}$ such that 
$
\blambda _0=\mathbb{LH}_{\mathcal{S}_{t_0}}(((u_0,t_0),\bxi_0),\mu_0)
$
and 
$1/{\mu_0}$ is one of the non-zero momentary nullcone
principal curvatures 
$\kappa_N(\mathcal{S}_{t_0})_i((u_0,t_0),\bxi_0), (i=1,\dots ,s).$
\par\noindent
{\rm (3)} Under the condition {\rm (2)}, ${\rm rank}\, {\mathcal H}(h_{(t_0,\sblambda _0)})(u_0)=0$ if and only if
$p_0=\bX(u_0,t_0)$ is a non-parabolic momentary $\bxi _0$-nullcone umbilical point.
\end{Pro}
\demo
(1) We denote that $p_0=\bX(u_0,t_0).$ The condition $h_{(t_0,\sblambda _0})(u_0)=\langle \bX(u_0,t_0),{\blambda
_0}\rangle+1 =0$
means that
\begin{eqnarray*}
\langle \bX(u_0,t_0)-\lambda _0,\bX(u_0.t_0)-\blambda_0\rangle &=&\langle\bX(u_0,t_0),\bX(u_0,t_0)\rangle-2\langle\bX(u_0,t_0),\blambda _0\rangle+\langle\blambda _0,\blambda_0\rangle \\
&=&-2(1+\langle\bX(u_0,t_0),\blambda _0\rangle)=0,
\end{eqnarray*}
so that
$\bX(u_0,t_0)-{\blambda _0}\in \Lambda^*.$
Since $\partial h_{(t_0,\sblambda _0)}/\partial u_i(u)=\langle \bX_{u_i}(u,t_0), {\blambda _0}\rangle $
and $\langle \bX_{u_i},\bX\rangle=0,$
we have $\langle \bX_{u_i}(u,t_0),\blambda _0\rangle=-\langle \bX_{u_i}(u,t_0)-\blambda _0\rangle$.
Therefore, $\partial h_{(t_0,\sblambda _0)}/\partial u_i(u_0)=0$
if and only if
$\bX (u_0,t_0)-{\blambda _0}\in N_{p_0}M.$
On the other hand, the condition $h_{(t_0,\sblambda _0)} (u_0)=\langle \bX(u_0,t_0),\blambda _0\rangle+1=0$
implies that $\langle\bX(u_0,t_0),\bX(u_0,t_0)-\blambda_0\rangle =0$.
This means that $\bX(u_0,t_0)-\blambda _0\in T_{p_0}AdS^{n+1}.$
Hence
$h_{(t_0,\sblambda _0)}(u_0)=\partial h_{(t_0,\sblambda _0)}/\partial u_i(u_0)=0$ $(i=1,\dots, s)$
if and only if
$\bX(u_0,t_0)-{\blambda _0}\in N_{p_0}(\mathcal{S}_{t_0})\cap \Lambda^*\cap T_{p_0}AdS^{n+1}.$
Then we denote that
$\bv=\bX (u_0,t_0)-{\blambda _0}\in N_{p_0}(\mathcal{S}_{t_0})\cap \Lambda^*\cap T_{p_0}AdS^{n+1}.$
If $\langle \bn^T(u_0,t_0),\bv\rangle =0,$ then $\bn^T(u_0,t_0)$ belongs to
a lightlike hyperplane in the Lorentz space $T_{p_0}AdS^{n+1},$ so that $\bn^T(u_0,t_0)$ is lightlike or spacelike.
This contradiction to the fact that $\bn^T(u_0,t_0)$ is a timelike unit vector. Thus,
$\langle \bn^T(u_0,t_0),\bv\rangle \not=0.$ 
We set
\[
\bxi_0=\frac{-1}{\langle \bn^T(u_0,t_0),\bv\rangle}\bv -\bn^T(u_0,t_0).
\]
Then we have
\begin{eqnarray*}
\langle \bxi_0,\bxi_0\rangle &=& -2\frac{-1}{\langle \bn^T(u_0,t_0),\bv\rangle} \langle \bn^T(u_0,t_0),\bv\rangle-1=1 \\
\langle \bxi_0,\bn^T(u_0,t_0)\rangle &=& \frac{-1}{\langle \bn^T(u_0,t_0),\bv\rangle} \langle \bn^T(u_0,t_0),\bv\rangle+1=0.
\end{eqnarray*}
This means that $\bxi_0\in N_1[\mathcal{S}_{t_0}]_{p_0}.$
Since $-\bv=\langle \bn^T(u_0,t_0),\bv\rangle(\bn^T(u_0,t_0)+\bxi_0),$
we have 
${\blambda _0}=\bX(u_0,t_0)+\mu_0\mathbb{NG}(\mathcal{S}_{t_0})((u_0,t_0)\bxi_0)$, where
$p_0=\bX(u_0,t_0)$ and $\mu_0=\langle \bn^T(u_0,t_0),\bv\rangle.$
For the converse assertion,  suppose that $\blambda_0=\bX(u_0,t_0)+\mu_0\mathbb{NG}(\mathcal{S}_{t_0})((u_0,t_0),\bxi_0).$
Then $\blambda_0-\bX(u_0,t_0)\in N_{p_0}(\mathcal{S}_{t_0}))\cap \Lambda ^*$ and
$\langle\blambda_0-\bX(u_0,t_0),\bX(u_0,t_0)\rangle=\langle \mu_0\mathbb{NG}(\mathcal{S}_{t_0})(p_0,\bxi_0),\bX(u_0)\rangle=0.$
Thus we have $\blambda_0-\bX(u_0)\in N_{p_0}(\mathcal{S}_{t_0})\cap \Lambda ^*\cap T_{p_0}AdS^{n+1}.$
By the previous arguments, these conditions are equivalent to the condition that
$h_{(t_0,\sblambda _0)}(u_0)=\partial h_{(t_0,\sblambda _0)}/\partial u_i(u_0)=0$ $(i=1,\dots, s)$.

\par
(2) By a straightforward calculation, we have
\[
\frac{\partial ^2 h_{(t_0,\sblambda_0)}}{\partial u_i\partial u_j}(u)
=\langle\bX _{u_iu_j}(u,t_0),\blambda _0\rangle.
\]
Under the conditions ${\blambda _0}=\bX(u_0)+\mu_0(\bn^T(u_0)+\bxi_0)$,
we have
\[
\frac{\partial ^2 h_{(t_0,\sblambda_0)}}{\partial u_i\partial u_j}(u_0)
=\langle \bX _{u_iu_j}(u_0,t_0),\bX(u_0,t_0)\rangle +\mu_0\langle\bX_{u_iu_j}(u_0,t_0), (\bn^T(u_0,t_0)+\bxi_0)\rangle .
\]
Since $\langle \bX_{u_i},\bX\rangle =0,$ we have $\langle\bX_{u_iu_j},\bX\rangle=-\langle \bX_{u_i},\bX_{u_j}\rangle.$
Therefore, we have
\[
\left(\frac{\partial ^2 h_{(t_0,\sblambda_0)}}{\partial u_i\partial u_\ell}(u_0)\right)\left(g^{j\ell}(u_0,t_0)\right)
=\left(\mu_0 h^j_i(\mathcal{S}_{t_0})((u_0,t_0),\bxi_0)-\delta ^j_i\right).
\]
Thus, ${\rm det}{\mathcal H}(h_{(t_o,\sbxi_0)})(u_0)=0$ if and only if 
$1/\mu_0$ is an eigenvalue of $(h^i_j(\mathcal{S}_{t_0})((u_0,t_0),\bxi_0)),$ which is equal to
one of the momentary nullcone principal curvatures $\kappa _N(\mathcal{S}_{t_0})_i((u_0,t_0),\bxi_0),$ $(i=1,\dots ,s)$.
\par
(3) By the above calculation,  ${\rm rank}\, {\mathcal H}(h_{(t_0,\sblambda _0)})(u_0)=0$ if and only if
$$
(h^i_j(\mathcal{S}_{t_0})((u_0,t_0),\bxi_0))=\frac{1}{\mu _0}(\delta ^j_i),
$$
where $1/\mu_0=\kappa _N(\mathcal{S}_{t_0})_i((u_0,t_0),\bxi_0),$ $(i=1,\dots ,s)$. This means that $p_0=\bX(u_0,t_0)$ is a non-parabolic momentary $\bxi_0$-nullcone umbilical point.
\enD

\section{Graph-like big fronts}
In this section we briefly review the theory of
graph-like Legendrian unfoldings.
Graph-like Legendrian unfoldings belong to a special class of big Legendrian 
submanifolds (for detail, see \cite{Izumiya93,Izumiya-Takahashi,Izumiya-Takahashi2,Izumiya-Takahashi3,Zakalyukin95}).
Recently there appeared a survey article \cite{Graph-like} on the theory of graph-like Legendrian unfoldings.
Let ${\mathcal F} :(\R^k\times (\R^m\times\R),0)\to (\R,0)$ be a function germ.
We say that ${\mathcal F}$ is a {\it graph-like Morse family of hypersurfaces}
if $
(\mathcal{F}, d_q\mathcal{F}):(\R^k\times(\R^m\times \R),0)\to (\R\times \R^k,0)$
is a non-singular and 
$(\partial {\mathcal F}/\partial t)(0)\not= 0,$
where
$$
d_q\mathcal{F}(q,x,t)=\left(\frac{\partial \mathcal{F}}{\partial q_1}(q,x,t), \dots ,
\frac{\partial \mathcal{F}}{\partial q_k}(q,x,t)\right).
$$
Moreover, we say that ${\mathcal F}$ is {\it non-degenerate} if
$(\mathcal{F}, d_q\mathcal{F})|_{\R^k\times(\R^m\times \{0\})}$ is non-singular.
For a graph-like Morse family of hypersurfaces $\mathcal{F},$ $\Sigma _*(\mathcal{F})=(\mathcal{F}, d_q\mathcal{F})^{-1}(0)$ is a
smooth $m$-dimensional submanifold germ of $(\R^k\times(\R^m\times \R),0).$
We now consider the space of $1$-jets $J^1(\R^m,\R)$ with the canonical coordinates $(x_1,\dots ,x_m,t,p_1,\dots ,p_m)$ such
that the canonical contact form is $\theta =dt-\sum_{i=1}^m p_idx_i.$
We define a mapping $\Pi:J^1(\R^m,\R)\lon T^*\R^m$ by
$\Pi(x,t,p)=(x,p),$ where $(x,t,p)=(x_1,\dots, x_m,t,p_1,\dots ,p_m).$ 
Here, $T^*\R^m$ is a symplectic manifold with the
canonical symplectic structure $\omega=\sum_{i=1}^m dp_i\wedge dx_i$ (cf. \cite{Arnold1}).
We define a mapping $\mathscr{L}_{\mathcal{F}}:(\Sigma_{*}(\mathcal{F}),0) \to J^1(\R^m,\R)$ by 
\[
\mathscr{L}_{\mathcal F}(q,x,t)=\left(x,t,-\frac{\displaystyle \frac{\partial \mathcal{F}}{\displaystyle\partial x_1}(q,x,t)}{\frac{\displaystyle \partial \mathcal{F}}{\displaystyle\partial t}(q,x,t)},\dots , -\frac{\displaystyle\frac{\partial \mathcal{F}}{\partial x_m}(q,x,t)}{\frac{\displaystyle\partial \mathcal{F}}{\displaystyle\partial t}(q,x,t)},\right).
\]
It is easy to show that $\mathscr{L}_{\mathcal F}(\Sigma _*({\mathcal F}))$ is a Legendrian submanifold germ (cf., \cite{Arnold1}), which is called a {\it graph-like Legendrian unfolding germ.} 
We call $\overline{\pi}|_{\mathscr{L}_{\mathcal F}(\Sigma _*({\mathcal F}))}:\mathscr{L}_{\mathcal F}(\Sigma _*({\mathcal F}))\lon \R^m\times \R$  a {\it graph-like Legendrian map germ}, where
$\overline{\pi}:J^1(\R^m,\R)\lon \R^m\times\R$ is the canonical projection.
We also call $W(\mathscr{L}_{\mathcal F}(\Sigma _*({\mathcal F})))=\overline{\pi}(\mathscr{L}_{\mathcal F}(\Sigma _*({\mathcal F})))$ a {\it graph-like big front}
of $\mathscr{L}_{\mathcal F}(\Sigma _*({\mathcal F})).$
We say that ${\mathcal F}$ is a {\it graph-like generating family} of
$\mathscr{L}_{\mathcal F}(\Sigma _*({\mathcal F})).$
Moreover, we call $W_t(\mathscr{L}_{\mathcal F}(\Sigma _*({\mathcal F})))=\pi _1(\pi _2^{-1}(t)\cap W(\mathscr{L}_{\mathcal F}(\Sigma _*({\mathcal F})))$
a {\it momentary front} for each $t\in (\R,0),$ where $\pi _1:\R^m\times \R\lon \R^m$ and $\pi _2:\R^m\times \R\lon \R$ are the canonical projections.
The {\it discriminant set of the family} $\{W_t(\mathscr{L}_{\mathcal F}(\Sigma _*({\mathcal F})))\}_{t\in (\R,0)}$ is defined by the union of the {\it caustic} 
$$C_{\mathscr{L}_{\mathcal F}(\Sigma _*({\mathcal F}))}=\pi _1( \Sigma (W (\mathscr{L}_{\mathcal F}(\Sigma _*({\mathcal F}))))$$  and the {\it Maxwell stratified set} 
$$M_{\mathscr{L}_{\mathcal F}(\Sigma _*({\mathcal F}))}=\pi _1(SI_{W(\mathscr{L}_{\mathcal F}(\Sigma _*({\mathcal F})))}),$$ where $\Sigma (W (\mathscr{L}_{\mathcal F}(\Sigma _*({\mathcal F})))$ is the critical value set of $\overline{\pi}|_{\mathscr{L}_{\mathcal F}(\Sigma _*({\mathcal F}))}$
and $SI_{W(\mathscr{L}_{\mathcal F}(\Sigma _*({\mathcal F})))}$ is the closure of the self intersection set of $W(\mathscr{L}_{\mathcal F}(\Sigma _*({\mathcal F}))).$
\par
 We now define equivalence relations among graph-like Legendrian unfoldings.
Let ${\mathcal F} :(\R^k\times (\R^m\times\R),0)\to (\R,0)$ and ${\mathcal G} :(\R^{k}\times (\R^m\times\R),0)\to (\R,0)$ be graph-like Morse families of hypersurfaces.
We say that $\mathscr{L}_{\mathcal F}(\Sigma _*({\mathcal F}))$ and $\mathscr{L}_{\mathcal G}(\Sigma _*({\mathcal G}))$
are {\it Legendrian equivalent} if there exist a diffeomorphism germ $\Phi:(\R^m\times \R,\overline{\pi}(p ))\lon (\R^m\times \R,\overline{\pi}(p' ))$ and
a contact diffeomorphism germ $\widehat{\Phi}:(J^1(\R^m,\R),p)\lon (J^1(\R^m,\R),p')$ such that $\overline{\pi}\circ\widehat{\Phi}=\Phi\circ \overline{\pi}$
and $\widehat{\Phi}(\mathscr{L}_{\mathcal F}(\Sigma _*({\mathcal F})))=(\mathscr{L}_{\mathcal G}(\Sigma _*({\mathcal G}))),$
where $p=\mathscr{L}_{\mathcal F}(0)$ and $p'=\mathscr{L}_{\mathcal G}(0).$
We also say that $\mathscr{L}_{\mathcal F}(\Sigma _*({\mathcal F}))$ and $\mathscr{L}_{\mathcal G}(\Sigma _*({\mathcal G}))$
are {\it $S.P^+$-Legendrian equivalent} if these are Legendrian equivalent by a diffeomorphism germ $\Phi:(\R^m\times \R,\overline{\pi}(p ))\lon (\R^m\times \R,\overline{\pi}(p' ))$ 
of the form $\Phi (x,t)=(\phi _1(x),t+\alpha (x))$ and 
a contact diffeomorphism germ $\widehat{\Phi}:(J^1(\R^m,\R),p)\lon (J^1(\R^m,\R),p')$ with $\overline{\pi}\circ\widehat{\Phi}=\Phi\circ \overline{\pi}.$
Moreover, graph-like big fronts $W(\mathscr{L}_{\mathcal F}(\Sigma _*({\mathcal F})))$ and $W(\mathscr{L}_{\mathcal G}(\Sigma _*({\mathcal G})))$
are {\it $S.P^+$-diffeomorphic} if there exists a diffeomorphism germ $\Phi:(\R^m\times \R,\overline{\pi}(p ))\lon (\R^m\times \R,\overline{\pi}(p' ))$ 
of the form $\Phi (x,t)=(\phi _1(x),t+\alpha (x))$ such that $\Phi(W(\mathscr{L}_{\mathcal F}(\Sigma _*({\mathcal F}))))=W(\mathscr{L}_{\mathcal G}(\Sigma _*({\mathcal G})))$
as set germs.
By definition, if $\mathscr{L}_{\mathcal F}(\Sigma _*({\mathcal F}))$ and $\mathscr{L}_{\mathcal G}(\Sigma _*({\mathcal G}))$
are $S.P^+$-Legendrian equivalent, then $W(\mathscr{L}_{\mathcal F}(\Sigma _*({\mathcal F})))$ and $W(\mathscr{L}_{\mathcal G}(\Sigma _*({\mathcal G})))$
are $S.P^+$-diffeomorphic. The converse assertion holds generically \cite{Graph-like,GeomLag14}.
\begin{Pro}[\cite{GeomLag14}] Suppose that the sets of critical points of $\overline{\pi}|_{\mathscr{L}_{\mathcal F}(\Sigma _*({\mathcal F}))}, \overline{\pi}|_{\mathscr{L}_{\mathcal G}(\Sigma _*({\mathcal G}))}$
 are nowhere dense respectively.
 Then $\mathscr{L}_{\mathcal F}(\Sigma _*({\mathcal F}))$ and $\mathscr{L}_{\mathcal G}(\Sigma _*({\mathcal G}))$are $S.P^+$-Legendrian equivalent if and only if 
$
 W(\mathscr{L}_{\mathcal F}(\Sigma _*({\mathcal F})))$ and $W(\mathscr{L}_{\mathcal G}(\Sigma _*({\mathcal G})))
$
 are $S.P^+$-diffeomorphic.
\end{Pro}
\par
We remark that if $W(\mathscr{L}_{\mathcal F}(\Sigma _*({\mathcal F})))$ and $W(\mathscr{L}_{\mathcal G}(\Sigma _*({\mathcal G})))$
 are $S.P^+$-diffeomorphic by  a diffeomorphism germ $\Phi:(\R^m\times \R,\overline{\pi}(p ))\lon (\R^m\times \R,\overline{\pi}(p' ))$,
 then $$\Phi ( C_{\mathscr{L}_{\mathcal F}(\Sigma _*({\mathcal F}))}\cup M_{\mathscr{L}_{\mathcal F}(\Sigma _*({\mathcal F}))})=
 C_{\mathscr{L}_{\mathcal G}(\Sigma _*({\mathcal G}))}\cup M_{\mathscr{L}_{\mathcal G}(\Sigma _*({\mathcal G}))}.
 $$
\par
For a graph-like Morse family of hypersurfaces $\mathcal{F}:(\R^k\times (\R^m\times\R),0)\to (\R,0),$ by the implicit function theorem, there exist function germs
$F:(\R^k\times \R^m,0)\to (\R,0)$ and $\lambda :(\R^k\times (\R^m\times \R),0)\lon \R$ with $\lambda (0)\not= 0$
such that ${\mathcal F}(q,x,t)
=\lambda(q,x,t)( F(q,x)-t).$
We have shown in \cite{Graph-like} that
$\mathcal{F}$ is a graph-like Morse family of hypersurfaces if and only if
$F$ is a Morse family of functions.
Here we say that $F:(\R^k\times\R^m,0)\lon (\R,0)$ is a {\it Morse family of functions} if
\[
dF_q=\left(\frac{\partial F}{\partial q_1},\dots, \frac{\partial F}{\partial q_k}\right):(\R^k\times\R^m,0)\lon \R^k
\]
is non-singular.
We consider a graph-like Morse family of hypersurfaces 
\[
\mathcal{F}(q,x,t)=\lambda (q,x,t)(F(q,x)-t).
\]
In this case 
$
\Sigma _*(\mathcal{F})=\{(q,x,F(q,x))\in (\R^k\times (\R^m\times\R),0)\ |\ (q,x)\in C(F)\},
$
where
\[
C(F)=\left\{ (q,x)\in (\R^k\times\R^m,0)\ \Bigm|\ \frac{\partial F}{\partial q_1}(q,x)=\cdots =\frac{\partial F}{\partial q_k}(q,x)=0\ \right\}.
\]
Moreover, we define a map germ
$
L(F):(C(F),0)\lon T^*\R^m
$
by 
\[
L(F)(q,x)=\left(x,\frac{\partial F}{\partial x_1}(q,x),\dots ,\frac{\partial F}{\partial x_m}(q,x)\right)
\]
It is known that $L(F)(C(F))$ is a Lagrangian submanifold germ (cf., \cite{Arnold1})
for the canonical symplectic structure.
In this case $F$ is said to be a {\it generating family} of the Lagrangian submanifold germ $L(F)(C(F)).$
We remark that $\Pi (\mathscr{L}_{\mathcal F}(\Sigma _*({\mathcal F})))=L(F)(C(F))$ and the graph-like big front
$W(\mathscr{L}_{\mathcal F}(\Sigma _*({\mathcal F})))$ is the graph of $F|C(F).$
Here we call $\pi|_{L(F)(C(F))}:L(F)(C(F))\lon \R^m$ a {\it Lagrangian map germ}, where
$\pi :T^*\R^m\lon \R^m$ is the canonical projection.
Then the set of critical values of $\pi|_{L(F)(C(F))}$ is called a {\it caustic} of $L(F)(C(F))=\Pi(\mathscr{L}_{\mathcal{F}}(\Sigma _*(\mathcal{F})))$ in the theory of Lagrangian singularities,
which is denoted by $C_{L(F)(C(F))}.$
By definition, we have  $C_{L(F)(C(F))}=C_{\mathscr{L}_{\mathcal F}(\Sigma _*({\mathcal F}))}.$
\par
Let ${\mathcal F}, \mathcal{G} :(\R^k\times (\R^m\times\R),0)\to (\R,0)$ be graph-like Morse families of hypersurfaces.
We say that $\Pi(\mathscr{L}_{\mathcal{F}}(\Sigma _*(\mathcal{F})))$ and $\Pi(\mathscr{L}_{\mathcal{G}}(\Sigma _*(\mathcal{G})))$
are {\it Lagrangian equivalent} if there exist a diffeomorphism germ $\Psi:(\R^m, \pi\circ\Pi(p ))\lon (\R^m,\pi\circ\Pi(p' ))$ and
a symplectic diffeomorphism germ $\widehat{\Psi}:(T^*\R^m,\Pi( p))\lon (T^*\R^m,\Pi( p'))$ such that $\pi\circ\widehat{\Psi}=\Psi\circ \pi$
and $\widehat{\Psi}(\Pi(\mathscr{L}_{\mathcal{F}}(\Sigma _*(\mathcal{F}))))=\Pi(\mathscr{L}_{\mathcal{G}}(\Sigma _*(\mathcal{G}))),$
where $p=\mathscr{L}_{\mathcal F}(0)$ and $p'=\mathscr{L}_{\mathcal G}(0).$
By definition, if $\Pi(\mathscr{L}_{\mathcal{F}}(\Sigma _*(\mathcal{F})))$ and $\Pi(\mathscr{L}_{\mathcal{G}}(\Sigma _*(\mathcal{G})))$ are  Lagrangian equivalent, then the caustics $C_{\mathscr{L}_{\mathcal F}(\Sigma _*({\mathcal F}))}$ and $C_{\mathscr{L}_{\mathcal G}(\Sigma _*({\mathcal G}))}$ are diffeomorphic as set germs.
The converse assertion, however, does not hold (cf. \cite{GeomLag14}).
Recently, we have shown the following theorem (cf. \cite{Izumiya-Takahashi2,Graph-like, GeomLag14})
\begin{Th} With the same notations as the above, $\Pi(\mathscr{L}_{\mathcal{F}}(\Sigma _*(\mathcal{F})))$ and $\Pi(\mathscr{L}_{\mathcal{G}}(\Sigma _*(\mathcal{G})))$
are  Lagrangian equivalent if and only if 
$\mathscr{L}_{\mathcal F}(\Sigma _*({\mathcal F}))$ and $\mathscr{L}_{\mathcal G}(\Sigma _*({\mathcal G}))$
are $S.P^+$-Legendrian equivalent.
\end{Th}
We have the following corollary of Proposition 5.1 and Theorem 5.2.
\begin{Co} Suppose that the sets of critical points of $\overline{\pi}|_{\mathscr{L}_{\mathcal F}(\Sigma _*({\mathcal F}))}, \overline{\pi}|_{\mathscr{L}_{\mathcal G}(\Sigma _*({\mathcal G}))}$
 are nowhere dense, respectively.
 Then $\Pi(\mathscr{L}_{\mathcal{F}}(\Sigma _*(\mathcal{F})))$ and $\Pi(\mathscr{L}_{\mathcal{G}}(\Sigma _*(\mathcal{G})))$
are  Lagrangian equivalent  if and only if $W(\mathscr{L}_{\mathcal F}(\Sigma _*({\mathcal F})))$ and $W(\mathscr{L}_{\mathcal G}(\Sigma _*({\mathcal G})))$
 are $S.P^+$-diffeomorphic.
\end{Co}
There are the notions of Lagrangian stability of Lagrangian submanifold germs and $S.P^+$-Legendrian stability of graph-like Legendrian unfolding germs, respectively.
Here we do not use the exact definitions of those notions of stability, so that we omit to give the definitions.
For detailed properties of such stabilities, see \cite{Arnold1, Graph-like}.
We have the following corollary of Theorem 5.2.
\begin{Co} The graph-like Legendrian unfolding $\mathscr{L}_{\mathcal F}(\Sigma _*({\mathcal F}))$ is $S.P^+$-Legendrian stable if and only if 
the corresponding Lagrangian submanifold $\Pi(\mathscr{L}_{\mathcal{F}}(\Sigma _*(\mathcal{F})))$ is Lagrangian stable.
\end{Co}
\par
 Let ${\mathcal F}:(\R^k\times (\R^m\times\R),0)\to (\R,0)$ be a graph-like Morse family of hypersurfaces. We define $\overline{f}:(\R^k\times \R,0)\lon (\R,0)$ by
$\overline{f}(q,t)={\mathcal F}(q,0,t).$
For graph-like Morse families of hypersurfaces ${\mathcal F} :(\R^k\times (\R^m\times\R),0)\to (\R,0)$ and ${\mathcal G} :(\R^{k}\times (\R^m\times\R),0)\to (\R,0)$,
we say that $\overline{f}$ and $\overline{g}$ are {\it $S.P$-$\mathcal{K}$-equivalent} if there exist a function germ $\nu :(\R^k\times\R,0)\lon \R$ with $\nu(0)\not= 0$ and
a diffeomorphism germ $\phi :(\R^k\times \R,0)\lon (\R^k\times \R,0)$ of the form $\phi (q,t)=(\phi_1(q,t),t)$ such that $\overline{f}(q,t)=\nu (q,t)\overline{g}(\phi (q,t)).$
Although we do not give the definition of $S.P^+$-Legendrian stability, we give a corresponding notion for graph-like Morse family of hypersurfaces.
We say that ${\mathcal F}$ is an {\it infinitesimally $S.P^+$-$\mathcal{K}$-versal unfolding} of $\overline{f}$ if 
$$
{\cal E}_{k+1}=\left\langle \frac{\partial \overline{f}}{\partial q_1},\dots, \frac{\partial \overline{f}}{\partial q_k},\overline{f} \right\rangle _{{\cal E}_{k+1}}+
\left\langle \frac{\partial \overline{f}}{\partial t} \right\rangle _{\R}+
\left\langle \frac{\partial \mathcal{F}}{\partial x_1}|_{\R^k\times\{0\}\times \R} ,\dots
,\frac{\partial \mathcal{F}}{\partial x_m}|_{\R^k\times\{0\}\times\R} \right\rangle _{\R},
$$
where ${\cal E}_{k+1}$ is the local $\R$-algebra of $C^\infty$-function germs $(\R^k\times\R,0)\lon \R.$
It is known the following theorem in \cite{Izu95, Zakalyukin95}.
\begin{Th}
The graph-like Legendrian unfolding
$\mathscr{L} _\mathcal{F}(\Sigma _*(\mathcal{F}))$ is $S.P^+$-Legendre stable if and only if $\mathcal{F}$ is an infinitesimally $S.P^+$-${\cal K}$-versal unfolding of $\overline{f}.$
\end{Th}

In \cite{Graph-like} we have shown the following theorem.
\begin{Th} Let ${\mathcal F},\mathcal{G} :(\R^k\times (\R^m\times\R),0)\to (\R,0)$  be graph-like Morse families of hypersurfaces such that ${\mathscr{L}_{\mathcal F}(\Sigma _*({\mathcal F}))}, {\mathscr{L}_{\mathcal G}(\Sigma _*({\mathcal G}))}$ are $S.P^+$-Legendrian stable. Then the following conditions are equivalent{\rm :}
\par\noindent
{\rm (1)} $\mathscr{L}_{\mathcal{F}}(\Sigma _*(\mathcal{F}))$
and $\mathscr{L}_{\mathcal{G}}(\Sigma _*(\mathcal{G}))$ are $S.P^+$-Legendrian equivalent,
\par\noindent
{\rm (2)} $\overline{f}$ and $\overline{g}$ are $S.P$-$\mathcal{K}$-equivalent,
\par\noindent
{\rm (3)} $\Pi(\mathscr{L}_{\mathcal{F}}(\Sigma _*(\mathcal{F})))$ and $\Pi(\mathscr{L}_{\mathcal{G}}(\Sigma _*(\mathcal{G})))$ are Lagrangian equivalent,
\par\noindent
{\rm (4)} $W(\mathscr{L}_{\mathcal{F}}(\Sigma _*(\mathcal{F})))$ and $W(\mathscr{L}_{\mathcal{G}}(\Sigma _*(\mathcal{G})))$
are $S.P^+$-diffeomorphic.
\end{Th}
\section{Unfolded lightlike hypersrufaces}
\par
Returning to our situation,
we have the following proposition.
\begin{Pro} 
Let $H$ be the AdS-height function on $W.$
For any $((u,t),\blambda )\in \Delta ^*H^{-1}(0),$ the germ of $H$ at
$(u,\blambda )$
is a non-degenerate graph-like Morse family of hypersurfaces. 
\end{Pro}
\demo
We denote that
$$
\bX(u,t)=(X_{-1}(u,t),X_0(u,t),X_1(u,t),\dots ,X_n(u,t))\ {\rm and}\
\blambda =(\lambda _{-1},\lambda _0,\lambda _1,\dots ,\lambda _n).
$$
We define an open subset 
$U_{-1}^+=\{\blambda \in AdS^{n+1}\ |\ \lambda _{-1}>0\ \}$.
For any $\blambda \in U_{-1}^+,$ we have
\[
\lambda _{-1}=\sqrt{1-\lambda _0^2+\lambda _1^2+\cdots \lambda _n^2}.
\]
Thus, we have a local coordinate of $AdS^{n+1}$ given by $(\lambda _0,\lambda _1,\dots ,\lambda _n)$ on
$U_{-1}^+.$
By definition, we have
$$
H(u,t,\blambda )=-X_{-1}(u,t)\sqrt{1-\lambda _0^2+\sum_{i=1}^n \lambda _i^2}-X_0(u,t)\lambda_0+
X_1(u,t)\lambda_1 +\cdots 
+X_n(u,t)\lambda _n.
$$
We now prove that the mapping $$
\Delta^*H|(U\times \{t\}\times U_{-1}^+)=\left(H, \frac{\partial H}{\partial u_1},\dots ,\frac{\partial H}{\partial
u_s}\right):U\times\{t\}\times U_{-1}^+\lon \R\times \R^s
$$
is non-singular at $(u,t,\blambda )\in \Delta ^*H^{-1}(0)\cap (U\times\{t\}\times U_{-1}^+).$
Indeed, the Jacobian matrix of $\Delta ^*H|(U\times \{t\}\times U_{-1}^+)$ is given by
\newfont{\bg}{cmr10 scaled\magstep5}
\newcommand{\bigA}{\smash{\lower1.0ex\hbox{\bg A}}}
\[
\left(
\begin{array}{ccccc}
 &
X_{-1}\displaystyle{\frac{\lambda _0}{\lambda _{-1}}}-X_0 & -X_{-1}\displaystyle{\frac{\lambda _1}{\lambda _{-1}}}+X_1 &  \cdots  &
-X_{-1}\displaystyle{\frac{\lambda _n}{\lambda _{-1}}}-X_n \\
\bigA & X_{-1u_1}\displaystyle{\frac{\lambda _0}{\lambda _{-1}}}-X_{0u_1} &
-X_{-1u_1}\displaystyle{\frac{\lambda _1}{\lambda _{-1}}}+X_{1u_1}& \cdots  &-X_{-1u_1}\displaystyle{\frac{\lambda _n}{\lambda _{-1}}}-X_{nu_1}\\
&\vdots & \vdots & \ddots & \vdots \\
& X_{-1u_s}\displaystyle{\frac{\lambda _0}{\lambda _{-1}}}-X_{0u_s} &
-X_{-1u_s}\displaystyle{\frac{\lambda _1}{\lambda _{-1}}}+X_{1u_s}& \cdots  &-X_{-1u_s}\displaystyle{\frac{\lambda _n}{\lambda _{-1}}}-X_{nu_s}
\end{array}
\right) ,
\]
where
\begin{eqnarray*}
\bigA=
\left(\!\!
\begin{array}{ccc}
\langle \bX_{u_1} ,\blambda\rangle & \!\! \cdots\!\! & \langle \bX_{u_s},\blambda \rangle \\
\langle \bX_{u_1u_1},\blambda\rangle & \!\! \cdots\!\!  &
\langle \bX_{u_1u_s},\blambda \rangle \\
\vdots & \!\!\ddots\!\! & \vdots \\
\langle \bX_{u_su_1},\blambda \rangle &\!\! \cdots\!\! &
\langle \bX_{u_su_s},\blambda\rangle 
\end{array}
\!\!
\right) .
\end{eqnarray*}
\newcommand{\bigB}{\smash{\lower1.0ex\hbox{\bg B}}}
We now show that
the rank of 
\[
\bigB=
\left(
\begin{array}{cccc}
X_{-1}\displaystyle{\frac{\lambda _0}{\lambda _{-1}}}-X_0 & -X_{-1}\displaystyle{\frac{\lambda _1}{\lambda _{-1}}}+X_1 &  \cdots  &
-X_{-1}\displaystyle{\frac{\lambda _n}{\lambda _{-1}}}-X_n \\
 X_{-1u_1}\displaystyle{\frac{\lambda _0}{\lambda _{-1}}}-X_{0u_1} &
-X_{-1u_1}\displaystyle{\frac{\lambda _1}{\lambda _{-1}}}+X_{1u_1}& \cdots  &-X_{-1u_1}\displaystyle{\frac{\lambda _n}{\lambda _{-1}}}-X_{nu_1}\\
\vdots & \vdots & \ddots & \vdots \\
 X_{-1u_s}\displaystyle{\frac{\lambda _0}{\lambda _{-1}}}-X_{0u_s} &
-X_{-1u_s}\displaystyle{\frac{\lambda _1}{\lambda _{-1}}}+X_{1u_s}& \cdots  &-X_{-1u_s}\displaystyle{\frac{\lambda _n}{\lambda _{-1}}}-X_{nu_s}
\end{array}
\right) 
\]
is $s+1$ at $(u,t,\blambda)\in \Sigma _*(H).$
Since $(u,t,\blambda)\in \Sigma _*(H),$ we have
\[
\blambda =\bX(u,t)+\mu\left(\bn^T(u,t)+\sum _{i=1}^{k-1}\xi _i\bn_i(u,t)\right)
\]
with $\sum_{i=1}^{k-1}\xi ^2_i=1,$ where 
$\{\bX,\bn^T,\bn^S_1,\dots ,\bn^S_{k-1}\}$ is a pseudo-orthonormal (local) frame of $N(M).$
Without the loss of generality, we assume that $\mu\not= 0$ and $\xi_{k-1}\not= 0.$
We denote that
\[
\bn^T(u,t)=^t\!\!(n^T_{-1}(u,t),n^T_0(u,t),\dots n^T_n(u,t)),\ 
\bn_i(u,t)=^t\!\!(n^i_{-1}(u,t),n_0^i(u,t),\dots n^i_n(u,t)).
\]
\newcommand{\bigC}{\smash{\lower1.0ex\hbox{\bg C}}}
It is enough to show that the rank of 
the matrix
\[
\bigC =
\left(
\begin{array}{cccc}
X_{-1}\displaystyle{\frac{\lambda _0}{\lambda _{-1}}}-X_0 & -X_{-1}\displaystyle{\frac{\lambda _1}{\lambda _{-1}}}+X_1 &  \cdots  &
-X_{-1}\displaystyle{\frac{\lambda _n}{\lambda _{-1}}}-X_n \\
 X_{-1u_1}\displaystyle{\frac{\lambda _0}{\lambda _{-1}}}-X_{0u_1} &
-X_{-1u_1}\displaystyle{\frac{\lambda _1}{\lambda _{-1}}}+X_{1u_1}& \cdots  &-X_{-1u_1}\displaystyle{\frac{\lambda _n}{\lambda _{-1}}}-X_{nu_1}\\
\vdots & \vdots & \ddots & \vdots \\
 X_{-1u_s}\displaystyle{\frac{\lambda _0}{\lambda _{-1}}}-X_{0u_s} &
-X_{-1u_s}\displaystyle{\frac{\lambda _1}{\lambda _{-1}}}+X_{1u_s}& \cdots  &-X_{-1u_s}\displaystyle{\frac{\lambda _n}{\lambda _{-1}}}-X_{nu_s} \\
n^T_{-1}\displaystyle{\frac{\lambda _0}{\lambda _{-1}}}-n^T_0 & -n^T_{-1}\displaystyle{\frac{\lambda _1}{\lambda _{-1}}}+n^T_1 &  \cdots  &
-n^T_{-1}\displaystyle{\frac{\lambda _n}{\lambda _{-1}}}-n^T_n \\
n^1_{-1}\displaystyle{\frac{\lambda _0}{\lambda _{-1}}}-n^1_0 & -n^1_{-1}\displaystyle{\frac{\lambda _1}{\lambda _{-1}}}+n^1_1 &  \cdots  &
-n^1_{-1}\displaystyle{\frac{\lambda _n}{\lambda _{-1}}}-n^1_n \\
\vdots & \vdots & \ddots & \vdots \\
n^{k-2}_{-1}\displaystyle{\frac{\lambda _0}{\lambda _{-1}}}-n^{k-2}_0 & -n^{k-2}_{-1}\displaystyle{\frac{\lambda _1}{\lambda _{-1}}}+n^{k-2}_1 &  \cdots  &
-n^{k-2}_{-1}\displaystyle{\frac{\lambda _n}{\lambda _{-1}}}-n^{k-2}_n 
\end{array}
\right) 
\]
is $n+1$ at $(u,t,\blambda)\in \Sigma _*(H).$
We denote that
\[
\ba_i=^t\!\!(x_i(u,t),x_{iu_1}(u,t),\dots x_{iu_s}(u,t),n^T_i(u,t),n^1_i(u,t),\dots ,n^{k-2}_i(u,t)).
\]
Then we have
\[
\bigC=\left(\ba_{-1}\frac{\lambda_0}{\lambda _{-1}}-\ba_0,-\ba_{-1}\frac{\lambda _1}{\lambda _{-1}}+\ba_1,\dots ,
-\ba _{-1}\frac{\lambda _n}{\lambda _{-1}}+\ba_n\right).
\]
It follows that
\begin{eqnarray*}
\det \bigC\!\!\!\!\!\!\!
&{}&=\frac{\lambda _{-1}}{\lambda _{-1}}\det (\ba_0,\ba_1,\dots,\ba_n)+\frac{\lambda_0}{\lambda_{-1}}\det (\ba_{-1}\ba_{1},\dots ,\ba_n)
\\
&{}&-\frac{\lambda _1}{\lambda _{-1}}(-1)\det(\ba_{-1},\ba_0,\ba_2,\dots ,\ba_n)-\cdots 
-\frac{\lambda_n}{\lambda _{-1}}(-1)^{n-1}\det(\ba _{-1}\ba_0,\ba_1,\dots ,\ba_{n-1}).
\end{eqnarray*}
Moreover, we define
$\delta _i=\det (\ba_{-1},\ba_0,\ba_1,\dots,\ba _{i-1},\ba_{i+1},\dots ,\ba_n)$ for $i=-1,0,1,\dots ,n$
and
$\ba=(-\delta _{-1},-\delta _0,-\delta _1,(-1)^2\delta _2,\dots ,(-1)^{n-1}\delta _n).$
Then we have
\[
\ba=\bX\wedge\bX_{u_1}\wedge\cdots \wedge \bX_{u_s}\wedge\bn^T\wedge\bn_1\wedge\cdots\wedge\bn_{k-2}.
\]
We remark that $\ba\not=0$ and $\ba=\pm\|\ba\|\bn_{k-1}.$
By the above calculation, we have
\begin{eqnarray*}
\det\bigC\!\!\!\!\!\!\!&{}&=\left\langle\left(\frac{\lambda{-1}}{\lambda_{-1}},\frac{\lambda_0}{\lambda _{-1}},\dots,\frac{\lambda_n}{\lambda _{-1}}\right), \ba\right\rangle =\frac{1}{\lambda_{-1}}\left\langle \bX(u)+\mu\left(\bn^T(u)+\sum_{i=1}^{k-1}\xi_i\bn_i(u)\right),\ba\right\rangle \\
&{}&=\frac{1}{\lambda _{-1}}\times\pm\mu\xi_{k-1}\|\ba\|=\pm \frac{\mu\xi_{k-1}\|\ba\|}{\lambda _{-1}}\not= 0.
\end{eqnarray*}
Therefore the Jacobi matrix of $\Delta^*H$
is non-singular at $(u,t,\blambda )\in \Delta ^*H^{-1}(0).$
\par
For other local coordinates of $AdS^{n+1}$, we can apply the same method for the proof as the above case.
Therefore, the AdS-height function $H$ is a non-degenerate big Morse family of hypersurfaces.
\par
On the other hand, we have
\[
\frac{\partial H}{\partial t}(u,t,\blambda)=\langle \bX_t(u,t),\blambda\rangle.
\]
Since $\bxi \in N^{AdS}_1[\mathcal{S}_t]_p=N_1^{AdS}(W)_p$ and $\bX_t(u,t)\in T_pW,$ we have
$\langle \bX_t(u,t) , \bxi\rangle =0.$
Moreover, we have $\langle \bX,\bX\rangle =-1,$ so that $\langle\bX_t(u,t),\bX(u,t)\rangle =0.$
Therefore, for $\blambda =\bX(u,t)+\mu(\bn^T(u,t)+\bxi),$ we have
\[
\frac{\partial H}{\partial t}(u,t,\blambda)=\langle \bX_t(u,t),\blambda\rangle
=\mu \langle \bX_t(u,t),\bn^T(u,t)\rangle.
\]
We remark that $\bn^T(u,t)$ is a timelike vector such that 
$\langle \bn^T(u,t),\bX_{u_i}(u,t)\rangle =0$, $(i=1,\dots s)$.
Since $\{\bX_t(u,t),\bX_{u_1}(u,t),\dots \bX_{u_s}(u,t)\}$
is a basis of the Lorentz space $T_pW$ and $\bn^T(u,t)\in T_pW,$
we have $\langle \bX_t(u,t),\bn^T(u,t)\rangle \not= 0.$
Moreover, $\blambda \notin W$ implies $\mu\not= 0.$
Thus we have $\partial H/\partial t(u,t)\not= 0$ for $\blambda =\bX(u,t)+\mu(\bn^T(u,t)+\bxi).$
This completes the proof.
\enD
\par
We also consider the local coordinate 
$U^+_{-1}$.
Since $H$ is a non-degenerate graph-like Morse family of hypersurfaces, we have a non-degenerate graph-like Legendrian unfolding
\[
\mathscr{L} _H:\Sigma _*(H)\lon J^1(U^+_{-1},I).
\]
By definition, we have
\[
\frac{\partial H}{\partial \lambda _0}((u,t),\blambda) \\
=X_{-1}(u)\displaystyle{\frac{\lambda _0}{\lambda _{-1}}}-X_0(u),\ \frac{\partial H}{\partial \lambda _i}((u,t),\blambda)=-X_{-1}(u)\displaystyle{\frac{\lambda _i}{\lambda _{-1}}}+X_i(u),
\]
$(i=1,\dots , n)$ and $\partial H/\partial t((u,t),\blambda)=\langle \bX_t(u,t),\blambda\rangle.$
It follows that
\begin{eqnarray*}
&{}&\left[\frac{\partial H}{\partial t}((u,t),\blambda):\frac{\partial H}{\partial \lambda _0}((u,t),\blambda):\frac{\partial H}{\partial \lambda _1}((u,t),\blambda):\cdots :\frac{\partial H}{\partial \lambda _n}((u,t),\blambda)\right]=
\\
&{}&[\langle \bX_t,\blambda\rangle:X_{-1}(u)\lambda _0-X_0(u)\lambda _{-1}:X_1(u)\lambda _{-1}-X_{-1}(u)\lambda _1:\cdots :X_n(u)\lambda _{-1}-X_{-1}(u)\lambda _n].
\end{eqnarray*}
We denote that 
\[
D_i(\bX,\blambda)=\det\begin{pmatrix}
X_{-1} & X_i \\
\lambda _{-1} & \lambda _i \\
\end{pmatrix},\
(i=0,1,\dots ,n).
\]
Then we have
\[
\mathscr{L}_H((u,t),\blambda )=\left(\blambda, t,-\frac{D_0((\bX,\blambda)}{\langle \bX_t,\blambda\rangle},\frac{D_1((\bX,\blambda)}{\langle \bX_t,\blambda\rangle},\dots ,\frac{D_n((\bX,\blambda)}{\langle \bX_t,\blambda\rangle}\right),
\]
where
\[
\Sigma _*(H)=\Bigl\{((u,t),\blambda)\ \Bigm|\ \blambda =\mathbb{LH}_{\mathcal{S}_t}(((u,t),\bxi),\mu)\ ((p,\bxi),\mu)\in 
N_1^{AdS}[\mathcal{S}_t]\times\R, p=\bX(u,t)\ \Bigr\}.
\]
We observe that $H$ is a graph-like generating family of the non-degenerate graph-like Legendrian unfolding 
$\mathscr{L}_H(\Sigma _*(H))$. 
Proposition 4.1 asserts that the graph-like big front $W(\mathscr{L}_H(\Sigma _*(H))$ of the non-degenerate graph-like Legendrian unfolding $\mathscr{L}_H(\Sigma _*(H))$ is given by
\[
\Bigl\{(\blambda,t)\in AdS^{n+1}\times I  \Bigm| \blambda =\mathbb{LH}_{\mathcal{S}_t}(((u,t),\bxi),\mu),\bxi\in N^{AdS}_1[\mathcal{S}_t]_p, p=\bX(u,t),\mu \in \R\ \Bigr\}.
\]
We define a mapping $\mathbb{LH}:N_1^{AdS}(W)\times \R\lon AdS^{n+1}\times I$ by
$$
\mathbb{LH}(\bX(u,t),\bxi,\mu)=(\mathbb{LH}_{\mathcal{S}_t}(\bX(u,t),\bxi,\mu),t),
$$
which is called an {\it unfolded lightlike hypersruface} of $W.$
We write $\mathbb{LH}_{(W,\mathcal{S})}=\mathbb{LH}(N_1^{AdS}(W)\times \R).$
Then we have $\mathbb{LH}_{(W,\mathcal{S})}=W(\mathscr{L}_H(\Sigma _*(H)),$
so that the image of the unfolded lightlike hypersruface of $W$ is the graph-like big front set of $\mathscr{L}_H(\Sigma _*(H)).$
Each momentary front is the lightlike hypersurface $\mathbb{LH}_{\mathcal{S}_t}(N_1^{AdS}[\mathcal{S}_t]\times\R)$, which is called a {\it momentary lightlike hypersruface} along the momentary space $\mathcal{S}_t.$ By assertion (2) of Proposition 4.1, a singular point of the momentary lightlike hypersruface  $\mathbb{LH}_{\mathcal{S}_t}(N_1^{AdS}[\mathcal{S}_t]\times\R)$ is a point $\blambda _0=\mathbb{LH}_{\mathcal{S}_{t_0}}(((u_0,t_0),\bxi_0,\mu_0)$
for 
$1/\mu _0 =\kappa  _N(\mathcal{S}_{t_0})_i((u_0,t_0),\bxi_0),$ $i=1,\dots ,s.$
Then we have the following corollary of Proposition 4.1.
\begin{Co}
 A singular point of $\mathbb{LH}_{(W,\mathcal{S})}$ is the point 
 $(\blambda, t)\in AdS^{n+1}\times I$ such that
$
\blambda =\mathbb{LH}_{\mathcal{S}_{t}}(((u,t),\bxi,\mu),$
where 
$1/\mu  =\kappa  _N(\mathcal{S}_{t})_i((u,t),\bxi),$ $i=1,\dots ,s.$
\end{Co}
\par
For a non-zero nullcone principal curvature $\kappa _N(\mathcal{S}_{t_0})_i((u_0,t_0),\bxi_0)\not= 0,$ we
have an open subset $O_i\subset N_1^{AdS}(W)$ such that $\kappa _N(\mathcal{S}_t)_i(\bX(u,t),\bxi)\not= 0$
for $(\bX(u,t),\bxi)\in O_i.$
Therefore, we have a non-zero nullcone principal curvature function $\kappa _N(\mathcal{S})_i:O_i\lon \R$.
We define a mapping
$
\mathbb{LF}_{\kappa_N(\mathcal{S}_t)_i} :O_i\cap N_1^{AdS}[\mathcal{S}_t]\lon AdS^{n+1}
$
by
\[
\mathbb{LF}_{\kappa_N(\mathcal{S}_t)_i}(\bX(u,t),\bxi)=\bX(u,t)+\frac{1}{\kappa_N(\mathcal{S}_t)_i(\bX(u,t),\bxi)}\mathbb{NG}((u,t),\bxi).
\]
We also  define
\[
\mathbb{LF}_{\mathcal{S}_t}=\bigcup_{i=1}^s \left\{\mathbb{LF}_{\kappa_N(\mathcal{S}_t)_i}(\bX(u,t),\bxi)\ |\ (\bX(u,t),\bxi)\in N_1^{AdS}[\mathcal{S}_t]\ \mbox{s.t.}\ \kappa _N(\mathcal{S}_t)_i(\bX(u,t),\bxi)\not= 0\right\} .
\]
We call $\mathbb{LF}_{\mathcal{S}_t}$ the {\it momentary lightlike focal set} along $\mathcal{S}_t=\bX(U\times\{t\})$ in $AdS^{n+1}.$
By definition, the momentary lightlike focal set along $\mathcal{S}_t=\bX(U\times\{t\})$ is the critical values set of the momentary lightlike hypersurface
$\mathbb{LH}_{\mathcal{S}_t}(N_1^{AdS}[\mathcal{S}_t]\times\R)$ along $\mathcal{S}_t$.
Moreover, an {\it unfolded lightcone focal set} of $(W,\mathcal{S})$ is defined to be
\[
\mathbb{LF}_{(W,\mathcal{S})}=\bigcup _{t\in I} \mathbb{LF}_{\mathcal{S}_{t}}\times \{t\} \subset AdS^{n+1}\times I.
\]
Then $\mathbb{LF}_{(W,\mathcal{S})}$ is the critical value set of $\mathbb{LH}$.

\section{Contact with lightcones}
\par
In this section we consider the geometric meanings of the singularities of momentary lightlike hypersrufaces in Anti-de Sitter space from the view point of the theory of contact of submanifolds with model hypersurfaces in \cite{mont1}.
We begin with the following basic observations.
\begin{Pro}
Let $\blambda _0\in AdS^{n+1}$ and $\mathcal{S}_{t_0}=\bX(U\times\{t_0\})$ a monetary space of $W=\bX(U\times I)$ without points
satisfying $K_N (\mathcal{S}_{t_0})(p,\bxi)= 0.$
Then $\mathcal{S}_{t_0}\subset
\Lambda ^{n+1}_{\lambda _0}\cap AdS^{n+1}$ 
if and only if $\blambda _0=\mathbb{LF}_{\mathcal{S}_{t_0}}$ is the momentary lightcone focal set.
In this case we have $\mathbb{LH}_{\mathcal{S}_{t_0}}(N^{AdS}_1[\mathcal{S}_{t_0}]\times \R)\subset \Lambda ^{n+1}_{\lambda _0}\cap AdS^{n+1}$
and $\mathcal{S}_{t_0}=\bX(U\times \{t_0\})$ is totally momentary nullcone umbilical.
\end{Pro}
\demo
By Proposition 3.1, $K_N(\mathcal{S}_{t_0})(p_0,\bxi_0)\not= 0$ if and only if
\[
\{(\bn^T+\bn^S), (\bn^T+\bn^S)_{u_1},
\dots , (\bn^T+\bn^S)_{u_{s}}\}
\]
is linearly independent for  $p_0=\bX(u_0,t_0)\in \mathcal{S}_{t_0}$ and $\bxi_0=\bn^S(u_0,t_0),$
where $\bn^S:\times I\lon N_1^{AdS}[\mathcal{S}_{t_0}]$ is a local section.
By the proof of the assertion (1) of Proposition 4.1, 
$\mathcal{S}_{t_0}\subset \Lambda ^{n+1}_{\lambda _0}\cap AdS^{n+1}$
if and only if $h_{\sblambda _0,t_0}(u)= 0$ for any $u\in U,$
where $h_{\sblambda _0,t_0}(u)=H(u,t_0,\blambda _0)$ is the AdS-height function on $\mathcal{S}_{t_0}.$
It also follows from  Proposition 4.1 that there exists a smooth function $\eta 
:U\times N_1^{AdS}[\mathcal{S}_{t_0}]\lon \R$ and section $\bn^S:U\times I\lon N_1^{AdS}[\mathcal{S}_{t_0}]$ such that
\[
\bX(u,t_0)=\blambda _0+\eta (u,\bn^S(u,t_0))(\bn^T(u,t_0)\pm\bn^S(u,t_0)).
\]
In fact, we have $\eta (u,\bn^S(u,t_0))=-1/\kappa _N(\mathcal{S}_{t_0})_i(p,\bxi)$ $i=1,\dots, s$, where
$p=\bX(u,t_0)$ and $\bxi=\bn^S(u,t_0).$
It follows that $\kappa _N(\mathcal{S}_{t_0})_i(p,\bxi)=\kappa _N(\mathcal{S}_{t_0})_j(p,\bxi),$ so that
$\mathcal{S}_{t_0}=\bX(U\times \{t_0\})$ is totally nullcone umbilical.
Therefore we have
\[
\mathbb{LH}_{\mathcal{S}_{t_0}}(u,\bn^S(u,t_0),\mu)=\blambda _0+(\mu +\eta (u,\bn^S(u,t_0))(\bn^T(u,t_0)\pm\bn^S(u,t_0)).
\]
Hence we have
$\mathbb{LH}_{\mathcal{S}_{t_0}}(N_1^{AdS}[\mathcal{S}_{t_0}]\times \R)\subset \Lambda ^{n+1}_{\lambda _0}\cap AdS^{n+1}.$
By definition, the critical value set of $\mathbb{LH}_{\mathcal{S}_{t_0}}(N_1^{AdS}[\mathcal{S}_{t_0}]\times \R)$ is the lightlike focal set
$\mathbb{LF}_{\mathcal{S}_{t_0}},$
which is equal to $\blambda _0$ by the previous arguments. 
\par
For the converse assertion, suppose that $\blambda _0=\mathbb{LF}_{\mathcal{S}_{t_0}}.$
Then we have 
\[
\blambda _0=\bX(u,t_0)+\frac{1}{\kappa _N(\mathcal{S}_{t_0})_i(\bX(u,t_0),\bxi)}\mathbb{NG}(\mathcal{S}_{t_0})(u,t_0,\bxi),
\]
for any $i=1,\dots ,s$ and $(p,\bxi)\in N_1^{AdS}[\mathcal{S}_{t_0}],$ where $p=\bX(u,t_0).$
Thus, we have 
\[
\kappa _N(\mathcal{S}_{t_0})_i(\bX(u,t_0),\bxi)=\kappa _N(\mathcal{S}_{t_0})_j(\bX(u,t_0),\bxi)
\]
for any $i,j=1,\dots ,s.$ This means that $\mathcal{S}_{t_0}$ is totally momentary nullcone umbilical.
Since $\mathbb{NG}(\mathcal{S}_{t_0})(u,t_0,\bxi)$ is null for any $(u,\bxi)$, we have $\bX(U\times\{t_0\})\subset \Lambda ^{n+1}_{\lambda_0}\cap AdS^{n+1}.$
This completes the proof.
\enD
\par
We now consider the relationship between the contact of a one parameter family of submanifolds with a submanifold and 
the $S.P$-${\mathcal K}$-classification of functions.  
Let $ U_i\subset \R^r$, ($i=1,2$) be open sets and $g_i:(U_i\times I, (\ou_i,t_i))\lon (\R^n,\bm{y}_i)$  immersion germs. We define $\overline{g}_i:(U_i\times I, (\ou_i,t_i))\lon (\R^n\times I,(\bm{y}_i,t_i))$
by $\overline{g}_i(\ou,t)=(g_i(\ou),t).$
We denote that $(\overline{Y}_i,(\bm{y}_i,t_i))=\overline{g}_i(U_i\times I),(\bm{y}_i,t_i)).$
Let 
$f_i:(\R^n,\bm{y}_i) \lon (\R,0)$ be submersion germs and denote that $(V(f_i),\bm{y}_i)=(f_i^{-1}(0),\bm{y}_i).$
We say that {\it the contact of $\overline{Y}_1$ with the trivial family of $V(f_1)$  
at $(\bm{y}_1,t_1)$} is of the {\it same type in the strict sense} as {\it the contact of 
$\overline{Y}_2$ with the trivial family of $V(f_2)$ at $(\bm{y}_2,t_2)$}
if there is a diffeomorphism germ
$\Phi:(\R^n\times I,(\bm{y}_1,t_1)) \lon (\R^n\times I,(\bm{y}_2,t_2))$ of the form $\Phi (\bm{y},t)=(\phi_1(\bm{y},t),t+(t_2-t_1))$
 such that $\Phi(\overline{Y}_1)=\overline{Y}_2$ and
$\Phi(V(f_1)\times I) = V(f_2)\times I$.
In this case we write 
$SK(\overline{Y}_1,V(f_1)\times I;(\bm{y}_1,t_1)) = SK(\overline{Y}_2,V(f_2)\times I;(\bm{y}_2,t_2))$.
We can show one of the parametric versions of Montaldi's theorem of contact between submanifolds as follows: 
\begin{Pro} 
We use the same notations as in the above paragraph. 
Then the following conditions are equivalent:
\par
{\rm (1)}
$
SK(\overline{Y}_1,V(f_1)\times I;(\bm{y}_1,t_1)) = SK(\overline{Y}_2,V(f_2)\times I;(\bm{y}_2,t_2))
$
\par
{\rm (2)}
$f_1 \circ g_1$ and $f_2 \circ g_2$ are $S.P$-${\mathcal K}$-equivalent {\rm (}i.e., there exists a diffeomorphism germ $\Psi :(U_1\times I,(\ou_1,t_1))\lon (U_2\times I,(\ou_2,t_2))$ of the form $\Psi (\ou,t)=(\psi_1 (\ou,t),t+(t_2-t_1))$ and a function germ $\lambda :(U_1\times I,(\ou_1,t_1))\lon \R$ with 
 $\lambda (\ou_1,t_1)\not= 0$ such that $(f_2\circ g_2)\circ \Phi (\ou,t) =\lambda (\ou,t)f_1\circ g_1(\ou,t)${\rm ).}
\end{Pro}
Since the proof of Proposition 7.2 is given by the arguments just along the line of the proof of the original theorem in \cite{mont1},
we omit the proof here.
\par
We now consider a function 
$
{\mathfrak h}_{\bm{\lambda}}:AdS^{n+1}\lon {\mathbb R}
$
defined by ${\mathfrak h}_{\bm{\lambda}}(\bx )=\langle \bx ,\bm{\lambda}\rangle +1,$
where $\bm{\lambda}\in AdS^{n+1}.$
For any $\bm{\lambda}_0\in AdS^{n+1}$, we have the Lorentzian tangent hyperplane $HP(\bm{\lambda}_0,-1)$
of de Sitter space $AdS^{n+1}$ at $\bm{\lambda}_0$, so that we have an AdS-lightcone 
\[
\mathfrak{h}_{\bm{\lambda}_0}^{-1}(0)=AdS^{n+1}\cap HP(\bm{\lambda}_0,-1)=LC^{AdS}(\bm{\lambda}_0 ).
\]
Moreover, we consider a point $\bm{\lambda}_0=\mathbb{LH}_{\mathcal{S}_{t_0}} (\bX(\ou_0,t_0),\bxi_0,\mu _0).$
Then we have 
$$
\mathfrak{h}_{\bm{\lambda}_0}\circ\bX (\ou_0,t_0)=H((u_0,t_0),\mathbb{LH}_{\mathcal{S}_{t_0}} (\bX(\ou_0,t_0),\bxi_0,\mu _0))=0.
$$
By Proposition 4.1, we also have relations that
$$
\frac{\partial \mathfrak{h}_{\bm{\lambda}_0}\circ\bX }{\partial u_i}(\ou_0,t_0)=\frac{\partial H}{\partial u_i}((\ou_0,t_0),\mathbb{LH}_{\mathcal{S}_{t_0}} (\bX(\ou_0,t_0),\bxi_0,\mu _0))=0.
$$
for $i=1,\dots ,s.$
This means that the AdS-lightcone $\mathfrak{h}_{\bm{\lambda}_0}^{-1}(0)=LC^{AdS}(\bm{\lambda}_0)$ is tangent to 
$\mathcal{S}_{t_0}=\bX (U\times \{t_0\})$ at $p_0=\bX (\ou_0,t_0).$
The AdS-lightcone $LC^{AdS}(\bm{\lambda}_0)$ is said to be a {\it tangent anti-de Sitter lightcone} (briefly, a {\it tangent
AdS-lightcone}) of $\mathcal{S}_{t_0}=\bX (U\times \{t_0\})$ at 
$p_0=\bX (\ou_0,t_0)$. We write that
$LC^{AdS}(\mathcal{S}_{t_0};p_0,\bxi_0,\mu _0)=LC^{AdS}(\bm{\lambda}_0),$ where $\bm{\lambda}_0=\mathbb{LH}_{\mathcal{S}_{t_0}} (\bX(\ou_0,t_0),\bxi_0,\mu _0).$
Then we have the following simple lemma.
\begin{Lem}
Let $\bX :U\times I\lon AdS^{n+1}$ be a world sheet in anti-de Sitter space. We consider two points $(p_1,\bxi_1,\mu_1),(p_2,\bxi_2,\mu_2)\in
N_1(\mathcal{S}_{t_0})\times \R ,$ where $p_i=\bX(\ou_i,t_0)$, $(i=1,2).$
Then  $\mathbb{LH}_{\mathcal{S}_{t_0}}(\bX(\ou_1,t_0),\bxi_1,\mu_1))=\mathbb{LH}_{\mathcal{S}_{t_0}} (\bX(\ou_2,t_0),\bxi_2,\mu_2))$ if and only if 
$$LC^{AdS}(\mathcal{S}_{t_0},p_1,\bxi_1,\mu_1)=LC^{AdS}(\mathcal{S}_{t_0},p_2,\bxi_2,\mu_2).$$
\end{Lem}
By the definition of unfolded lightlike hypersruface, 
\[
\mathbb{LH}(\bX(\overline{u}_1,t_1),\bxi_1,\mu_1)=\mathbb{LH}(\bX(\overline{u}_2,t_2),\bxi_2,\mu_2)
\] if and only if $t_1=t_2$ and $\mathbb{LH}_{\mathcal{S}_{t_1}} (\bX(\ou_1,t_1),\bxi_1,\mu_1)=\mathbb{LH}_{\mathcal{S}_{t_1}} (\bX(\ou_2,t_1),\bxi_2,\mu_2)$.
Eventually, we have tools for the study of the contact between world sheets and anti-de Sitter lightcones.
Since we have $h_{\bm{\lambda}}(\ou,t)=\mathfrak{h}_{\bm{\lambda}}\circ\bX(\ou,t),$ we have the following proposition as a
corollary of Proposition 7.2.
\begin{Pro}
Let $\bX_i : (U\times I,(\ou_i,t_i)) \lon (AdS^{n+1},p_i)$ $(i=1,2)$ be  world sheet germs with $W_i=\bX_i(U\times I)$ and $\bm{\lambda}_i=\mathbb{LH}_{\mathcal{S}_{t_i}} (\bX(\ou_i,t_i),\bxi_i,\mu _i).$
Then the following conditions are equivalent:
\par\noindent
{\rm (1)} $SK(\overline{W}_1, LC^{AdS}(\mathcal{S}_{t_1},p_1,\bxi_1,\mu_1)\times I;(p_1,t_1))=
SK(\overline{W}_2,LC^{AdS}(\mathcal{S}_{t_2},p_2,\bxi_2,\mu_2)\times I;(p_2,t_2)),$ 
\par\noindent
{\rm (2)} $h_{1,\bm{\lambda}_1}$ and $h_{2,\bm{\lambda}_2}$ are $S.P$-$\mathcal{K}$-equivalent.\\
\end{Pro}

\section{Caustics and Maxwell sets of world sheets}

In this section we apply the theory of graph-like Legendrian unfoldings to investigate the singularities of
the caustics and the Maxwell sets of world sheets.
 In \cite{Bousso, Bousso-Randall} Bousso and Randall gave an idea of caustics of world sheets in order to define the notion of holographic domains.
 The family of lightlike hypersrufaces $\{\mathbb{LH}_{\mathcal{S}_{t}}(N^{AdS}_1[\mathcal{S}_t]\times\R)\}_{t\in J}$ sweeps out a region in $AdS^{n+1}.$
A {\it caustic} of a world sheet is the union of the sets of critical values of lightlike hypersrufaces along momentary spaces $\{\mathcal{S}_t\}_{t\in I}.$
A {\it holographic domain} of the world sheet is the region where the light-sheets sweep out until {\it caustics}.
So this means that the boundary of the holographic domain consists the caustic of the world sheet.
The set of critical values of the lightlike hypersruface of a momentary space is the lightlike focal set of the momentary space.
Therefore the notion of caustics in the sense of Bousso-Randall is formulated as follows:
A {\it caustic of a world sheet} $(W,\mathcal{S})$ 
 is defined to be
\[
\displaystyle{C(W,\mathcal{S})=\bigcup _{t\in I} \mathbb{LF}_{\mathcal{S}_{t}}}=\pi_1(\mathbb{LF}_{(W,\mathcal{S})}),
\]
where $\pi_1:AdS^{n+1}\times I\lon AdS^{n+1}$ is the canonical projection.
We call $C(W,\mathcal{S})$ a {\it BR-caustic} of $(W,\mathcal{S}).$
By definition, we have $\Sigma (W(\mathscr{L}_{H}(\Sigma _*(H)))=\mathbb{LF} _{(W,\mathcal{S})},$ so that we have the following proposition.
\begin{Pro} Let $(W,\mathcal{S})$ be a world sheet in $AdS^{n+1}$ and $H:U\times I\times (AdS^{n+1}\setminus W)\lon \R$ 
the $AdS$-height function on $W.$ Then we have
$C(W,\mathcal{S})=C_{\mathscr{L}_{H}(\Sigma _*(H))}.$
\end{Pro}
\par
In \cite{Bousso,Bousso-Randall} the authors did not consider the Maxwell set of a world sheet.
However, the notion of Maxwell sets plays an important role in the cosmology which has been called a {\it crease set} by Penrose (cf. \cite{Penrose,Siino}).
Actually, the topological shape of the event horizon is determined by the crease set of lightlike hypersrufaces.
Here, we write $M(W,S)=M_{\mathscr{L}_H(\Sigma _*(H))}$ and call it a {\it BR-Maxwell set} of the world sheet $(W,\mathcal{S}).$ 
\par
Let $\bX_i:(U\times I, (\ou_i,t_i))\lon (AdS^{n+1},p_i)$, $(i=1,2)$ be germs of timelike embeddings such that
$(W_i,\mathcal{S}_{i})$ are world sheet germs, where $W_i=\bX_i(U\times I).$ 
For $\bm{\lambda}_i=\mathbb{LH}_{\mathcal{S}_{t_i}} (\bX(\ou_i,t_i),\bxi_i,\mu _i),$ let $H_i :(U\times I\times (AdS^{n+1}\setminus W_i) ,(\ou_i,t_i,\bm{\lambda}_i))\lon \R$
be $AdS$-height function germs.
We also write $h_{i,\bm{\lambda}_i}(\ou,t)=H_i(\ou,t,\bm{\lambda}_i).$
Since $$W(\mathscr{L}_{H_i}(\Sigma _*(H_i)))=\mathbb{LH}_{(W_i,\mathcal{S}_i)},$$
we can apply Theorem 5.2 and Corollary 5.3 to our case. Then we have the following theorem.

\begin{Th} Suppose that the set of critical points of $\overline{\pi}|_{\mathscr{L}_{H_i}(\Sigma _*(H_i))}$ are nowhere dense for $i=1,2$, respectively.
Then the following conditions are equivalent{\rm :}
\par\noindent
{\rm (1)} $(\mathbb{LH}_{(W_1,\mathcal{S}_1)}, \bm{\lambda}_1)$ and $(\mathbb{LH}_{(W_2,\mathcal{S}_2)}, \bm{\lambda}_2)$ are $S.P^+$-diffeomorphic,
\par\noindent
{\rm (2)} $\mathscr{L}_{H_1}(\Sigma _*(H_1))$ and $\mathscr{L}_{H_2}(\Sigma _*(H_2))$ are $S.P^+$-Legendrian equivalent,
\par\noindent
{\rm (3)} $\Pi(\mathscr{L}_{H_1}(\Sigma _*(H_1)))$ and $\Pi(\mathscr{L}_{H_2}(\Sigma _*(H_2))$ are Lagrangian equivalent.
\end{Th}
\par
We remark that conditions (2) and (3) are equivalent without any assumptions (cf. Theorem 5.2).
Moreover, if we assume that $\mathscr{L}_{H_i}(\Sigma _*(H_i))$ are $S.P^+$-Legendrian stable, then we can apply Proposition 7.4 and Theorem 5.6 to show the following theorem.
\begin{Th} Suppose that $\mathscr{L}_{H_i}(\Sigma _*(H_i))$ are $S.P^+$-Legendrian stable for $i=1,2,$ respectively.
Then the following conditions are equivalent{\rm :}
\par\noindent
{\rm (1)} $(\mathbb{LH}_{(W_1,\mathcal{S}_1)}, \bm{\lambda}_1)$ and $(\mathbb{LH}_{(W_2,\mathcal{S}_2)}, \bm{\lambda}_2)$ are $S.P^+$-diffeomorphic,
\par\noindent
{\rm (2)} $\mathscr{L}_{H_1}(\Sigma _*(H_1))$ and $\mathscr{L}_{H_2}(\Sigma _*(H_2))$ are $S.P^+$-Legendrian equivalent,
\par\noindent
{\rm (3)} $\Pi(\mathscr{L}_{H_1}(\Sigma _*(H_1)))$ and $\Pi(\mathscr{L}_{H_2}(\Sigma _*(H_2))$ are Lagrangian equivalent,
\par\noindent
{\rm (4)} $h_{1,\bm{\lambda}_1}$ and $h_{2,\bm{\lambda}_2}$ are $S.P$-$\mathcal{K}$-equivalent,
\par\noindent
{\rm (5)} $SK(\overline{W}_1, LC^{AdS}(\mathcal{S}_{t_1},p_1,\bm{\xi}_1,\mu_1)\times I;(p_1,t_1))=
SK(\overline{W}_2, LC^{AdS}(\mathcal{S}_{t_2},p_2,\bm{\xi}_2,\mu_2)\times I;(p_2,t_2)).$ 
\end{Th}
\par
By definition and Proposition 8.1, we have the following proposition.
\begin{Pro} If $\Pi(\mathscr{L}_{H_1}(\Sigma _*(H_1)))$ and $\Pi(\mathscr{L}_{H_2}(\Sigma _*(H_2))$ are Lagrangian equivalent,
then BR-caustics $C(W_1,\mathcal{S}_1)$, $C(W_2,\mathcal{S}_2)$ and BR-Maxwell sets $M(W_1,\mathcal{S}_1)$, $M(W_2,\mathcal{S}_2)$ are diffeomorphic as set germs, respectively.
\end{Pro}

\section{World hyper-sheets in $AdS^{n+1}$}
In this section we consider the case when $k=2.$
For an open subset $U\subset \R^n,$ let $\bX:U\times I\lon AdS^{n+1}$ be a timelike embedding such that $(W, \mathcal{S})$ is a world sheet.
In this case $(W,\mathcal{S})$ is said to be a {\it world hyper-sheet} in $AdS^{n+1}.$
Since the pseudo normal space $N_p(W)$ is a Lorentz plane, $N_p^{AdS}(W)$ is a spacelike line, so that $N_1^{AdS}(W)_p$ comprises two points.
For any $\bm{\xi}\in N_1^{AdS}(W)_p,$ we have $-\bm{\xi}\in N_1^{AdS}(W)_p.$
We define a pseudo normal section $\bm{n}^S(\ou,t)\in N_1^{AdS}(W)_p$ for $p=\bm{X}(\ou,t)$
by
\[
\bm{n}^S(\ou,t)=\frac{\bX(\ou,t)\wedge \bX_{u_1}(\ou,t)\wedge \dots \wedge \bX_{u_{n-1}}(\ou,t)\wedge \bX_t(u,t)}{\|\bX(\ou,t)\wedge \bX_{u_1}(\ou,t)\wedge \dots \wedge \bX_{u_{n-1}}(\ou,t)\wedge \bX_t(\ou,t)\|}.
\]
Therefore the momentary nullcone Gauss images 
\[
\mathbb{NG}(\mathcal{S}_{t_0},\pm\bm{n}^S):U\lon \Lambda ^*
\]
are given by $\mathbb{NG}(\mathcal{S}_{t_0},\pm\bm{n}^S)(\ou)=\bm{n}^T(\ou,t_0)\pm \bm{n}^S(\ou,t_0).$
Therefore we have the momentary nullcone shape operators
\[
S_N^\pm (\mathcal{S}_{t_0})_p=S_p(\mathcal{S}_{t_0};\pm\bm{n}^S)=-\pi ^t\circ d_p\mathbb{NG}(\mathcal{S}_{t_0},\pm\bm{n}^S):T_p\mathcal{S}_{t_0}\lon T_p\mathcal{S}_{t_0}.
\]
It follows that we have momentary nullcone principal curvatures 
$$\kappa _N^\pm(\mathcal{S}_{t_0})_i( p)=\kappa _N(\mathcal{S}_{t_0})(p,\pm \bm{n}^S(\ou,t_0)),\ ( i=1,\dots ,n-1).$$
Then the momentary lightlike hypersrufaces $\mathbb{LH}^\pm_{S_t}:U\times \R\lon AdS^{n+1}$ are given by
\[
\mathbb{LH}^\pm_{S_t}(\ou,\mu)=\bX(\ou,t)+\mu(\bm{n}^T(\ou,t)\pm\bm{n}^S(\ou,t))=\bX(\ou,t)+\mu\mathbb{NG}(\mathcal{S}_{t},\pm\bm{n}^S)(\ou).
\]
Moreover, the unfolded lightlike hypersrufaces $\mathbb{LH}^\pm:U\times \R\lon AdS^{n+1}\times I$ are given by
\[
\mathbb{LH}^\pm(\ou,\mu)=(\mathbb{LH}^\pm_{S_t}(\ou,\mu),t)=(\bX(\ou,t)+\mu\mathbb{NG}(\mathcal{S}_{t},\pm\bm{n}^S)(\ou),t).
\]
For the $AdS$-height function $H:U\times I\times AdS^{n+1}\lon \R$ on $(W,\mathcal{S}),$
$
\Sigma _*(H)=\Sigma ^+_*(H)\cup \Sigma ^-_*(H),
$ 
where
\[
\Sigma ^\pm_*(H)=\{((\ou,t),\bm{\lambda})\ |\ \bm{\lambda}=\mathbb{LH}^\pm_{\mathcal{S}_t}(\ou,t,\mu), \mu\in \R \}.
\]
Then the image of unfolded lightlike hypersrufaces is 
\[
\mathbb{LH}_W=\mathbb{LH}^+(U\times \R)\cup\mathbb{LH}^-(U\times \R)=W(\mathscr{L}_H(\Sigma_*(H))),
\]
which is the graph-like big front set of $\mathscr{L}_H(\Sigma_*(H)).$
The momentary lightlike focal sets along $\mathcal{S}_t$ are
\[
\mathbb{LF}^\pm_{\mathcal{S}_t}=\bigcup_{i=1}^{n-1} \left\{\mathbb{LF}^\pm_{\kappa^\pm_N(\mathcal{S}_t)_i}(\ou,t)\bigm |
(\ou,t)\in U\times I\ s.t.\ \kappa^\pm_N(\mathcal{S}_t)_i(\bX(\ou,t))\not= 0\right\},
\]
where
\[
\mathbb{LF}^\pm_{\kappa^\pm_N(\mathcal{S}_t)_i}(\ou,t)=\bX(\ou,t)+\frac{1}{\kappa^\pm_N(\mathcal{S}_t)_i(\bX(\ou,t))}\mathbb{NG}(\mathcal{S}_{t_0},\pm\bm{n}^S)(\ou).
\]
The unfolded lightcone focal set is
\[
\mathbb{LF}_{(W,\mathcal{S})}=\bigcup_{t\in I} \mathbb{LF}^+_{\mathcal{S}_t}\times \{t\}\cup \bigcup_{t\in I} \mathbb{LF}^-_{\mathcal{S}_t}\times \{t\}
\subset AdS^{n+1}\times I.
\]
In this case the BR-caustic is 
\[
C(W,\mathcal{S})=\pi_1(\mathbb{LF}_{(W,\mathcal{S})})=\bigcup_{t\in I}\mathbb{LF}^+_{\mathcal{S}_t}\cup \bigcup_{t\in I}\mathbb{LF}^-_{\mathcal{S}_t}.
\]
Moreover, the BR-Maxwell set is 
\[
M(W,\mathcal{S})=M_{\mathscr{L}_H(\Sigma_*(H))}=M_{\mathscr{L}_H(\Sigma^+_*(H))}\cup M_{\mathscr{L}_H(\Sigma^-_*(H))}.
\]

\section{World sheets in $AdS^3$}
In this section we consider world sheets in the $3$-dimensional anti de Sitter space as an example.
Let $(W, \mathcal{S})$ be a world sheet in $AdS^3$, which is parameterized by
a timelike embedding $\bm{\Gamma} :J\times I\lon AdS^{3}$ such that $\mathcal{S}_t=\bm{\Gamma} (J\times \{t\})$ for $t\in I.$
In this case we call $\mathcal{S}_t$ a {\it momentary curve}.
We assume that $s\in J$ is the arc-length parameter. Then  
$\bm{t}(s,t)=\bm{\gamma} _t'(s)$ is the unit spacelike tangent vector of $\mathcal{S} _t$, where $\bm{\gamma}_t(s)=\bm{\Gamma}(s,t).$
We have the unit pseudo-normal vector field $\bm{n}(s,t)$ of $W$ in $AdS^3$ defined by
\[
\bm{n}(s,t)=\frac{\bm{\Gamma}(s,t)\wedge \bm{t}(s,t)\wedge \bm{\Gamma}_t(s,t)}{\|\bm{\Gamma}(s,t)\wedge \bm{t}(s,t)\wedge \bm{\Gamma}_t(s,t)\|}.
\]
The unit timelike normal vector of $\mathcal{S} _t$ in $TW$ is defined to be $\bm{b}(s,t)=\bm{\Gamma}(s,t)\wedge \bm{n} (s,t)\wedge \bm{t}(s,t).$
We choose the orientation of $\mathcal{S} _t$ such that $\bm{b}(s,t)$ is adopted (i.e. ${\rm det}\, (\bm{\Gamma}(s,t), \bm{b}(s,t), \bm{e}_1,\bm{e}_2)> 0$).
Therefore,  
$\{\bm{\Gamma}(s,t),\bm{b}(s,t),\bm{n}(s,t),\bm{t}(s,t)\}$ is a {pseudo-orthonormal frame} along $W. $ 
On this moving frame, we can show the following
{\it Frenet-Serret type formulae} for $S_t$:
\[
	\setlength\arraycolsep{2pt}
	\left\{
	\begin{array}{ccl}
	\displaystyle{\frac{\partial \bm{\Gamma}}{\partial s}(s,t)} &=& \bm{t}(s,t), \\
	\displaystyle{\frac{\partial \bm{b}}{\partial s}(s,t)} &=& \tau_g(s,t) \bm{n}(s,t)-\kappa_g(s,t) \bm{t}(s,t), \\
    \displaystyle{\frac{\partial \bm{n}}{\partial s}(s,t)} &=& \tau_g(s,t) \bm{b}(s,t)-\kappa_n(s,t) \bm{t}(s,t) , \\
	\displaystyle{\frac{\partial \bm{t}}{\partial s}(s,t)} &=& \bm{\Gamma}(s,t)-\kappa_g(s,t) \bm{b}(s,t) + \kappa_n(s,t) \bm{n}(s,t),
	\end{array}
	\right. 
\]
where $\kappa_g(s,t) = \langle \frac{\partial \bm{t}}{\partial s}(s,t), \bm{b}(s,t) \rangle,$
	$\kappa_n(s,t) = \langle \frac{\partial \bm{t}}{\partial s}(s,t), \bm{n}(s,t) \rangle,$
	$\tau _g(s,t)= \langle \frac{\partial \bm{b}}{\partial s}(s,t),\bm{n}(s,t)\rangle.$
We call $\kappa _g(s,t)$ a {\it geodesic curvature}, $\kappa _n(s,t)$ a {\it normal curvature} and $\tau _g(s,t)$ a {\it geodesic torsion} of $\mathcal{S}_t$ respectively.
Then $\bm{b}(s,t_0)\pm \bm{n}(s,t_0)$ are lightlike. We have the momentary lightlike hypersrufaces
$\mathbb{LS}^\pm_{\mathcal{S}_{t_0}}:J\times\{t_0\}\times \R \lon AdS^3$ along $\mathcal{S}_{t_0}$ defined by
$
\mathbb{LS}^\pm_{\mathcal{S}_{t_0}}((s,t_0),u)=\bm{\Gamma} (s,t_0)+u(\bm{b}(s,t_0)\pm \bm{n}(s,t_0)).
$
Here, we use the notation $\mathbb{LS}^\pm_{\mathcal{S}_{t_0}}$ instead of $\mathbb{LH}^\pm_{\mathcal{S}_{t_0}}$
because the images of these mappings are lightlike surfaces.
We adopt $\bm{n}^T=\bm{b}$ and $\bm{n}^S=\bm{n}.$ By the Frenet-Serret type formulae, we have
\[
\frac{\partial (\bm{n}^T\pm \bm{n}^S)}{\partial s}(s,t)=\frac{\partial (\bm{b}\pm \bm{n})}{\partial s}(s,t)=\tau_g(s,t)(\bm{n}\pm\bm{b})(s,t)-(\kappa _g(s,t)\pm \kappa _n(s,t))\bm{t}(s,t).
\]
Therefore, we have $\kappa ^\pm (\mathcal{S}_t)(s,t)=\kappa _g(s,t)\pm \kappa _n(s,t).$
It follows that 
\[
\mathbb{LF}^\pm _{\mathcal{S}_{t_0}}=\left\{\bm{\Gamma}(s,t_0)+\frac{1}{\kappa _g(s,t_0)\pm \kappa _n(s,t_0)}(\bm{b}\pm \bm{t})(s,t_0)\bigm | s\in J, \kappa _g(s,t_0)\pm \kappa _n(s,t_0)\not= 0\right\}.
\]
We consider the $AdS$-height function $H:J\times I\times AdS^3\lon \R$.
Then we have
\begin{eqnarray*}
&{}& \frac{\partial H}{\partial s}(s,t,\bm{\lambda})=\langle \bm{t}(s,t),\bm{\lambda}\rangle, \\
&{}& \frac{\partial^2 H}{\partial s^2}(s,t,\bm{\lambda})=\langle (\bm{\Gamma}-\kappa _g\bm{b}+\kappa _n\bm{n})(s,t),\bm{\lambda}\rangle, \\
&{}& \frac{\partial^3 H}{\partial s^3}(s,t,\bm{\lambda})=\langle ((1+\kappa _g^2+\kappa _n^2)\bm{t}+(\kappa _n\tau _g-\kappa _g')\bm{b}+(\kappa _n'-\kappa _g\tau_g)\bm{n})(s,t),\bm{\lambda}\rangle.
\end{eqnarray*}
It follows that the following proposition holds. We write $H_{t_0}(s,\bm{\lambda})=H(s,t_0,\bm{\lambda}).$
\begin{Pro}
{\rm (1)}
$H_{t_0}(s,\bm{\lambda})=\partial H_{t_0}/\partial s (s,\bm{\lambda})=0$ if and only if there exists $ u\in \R$ such that  
 $\bm{\lambda}= \bm{\Gamma} (s,t_0)+u(\bm{b}(s,t_0)\pm \bm{n}(s,t_0))$
 \par\noindent
 {\rm (2)} $H_{t_0}(s,\bm{\lambda})=\partial H_{t_0}/\partial s (s,\bm{\lambda})=\partial^2 H_{t_0}/\partial s^2 (s,\bm{\lambda})=0$ if and only if $\kappa _g(s,t_0)\pm \kappa _n(s,t_0)\not= 0$ and \[
 \bm{\lambda}=\bm{\Gamma} (s,t_0)+\displaystyle{\frac{1}{\kappa _g(s,t_0)\pm \kappa _n(s,t_0)}}(\bm{b}(s,t_0)\pm \bm{n}(s,t_0)).
 \]
 \par\noindent
 {\rm (3)} $H_{t_0}(s,\bm{\lambda})=\partial H_{t_0}/\partial s (s,\bm{\lambda})=\partial^2 H_{t_0}/\partial s^2 (s,\bm{\lambda})=\partial^3 H_{t_0}/\partial s^3 (s,\bm{\lambda})=0$ if and only if $\kappa _g(s,t_0)\pm \kappa _n(s,t_0)\not= 0$, $((\kappa _n\pm\kappa_g)\tau _g\mp (\kappa _n'\pm\kappa _g'))(s_0,t_0)=0$  and 
 \[
 \bm{\lambda}=\bm{\Gamma} (s,t_0)+\displaystyle{\frac{1}{\kappa _g(s,t_0)\pm \kappa _n(s,t_0)}}(\bm{b}(s,t_0)\pm \bm{n}(s,t_0)).\]
 \par\noindent
 {\rm (4)} $H_{t_0}(s,\bm{\lambda})=\partial H_{t_0}/\partial s (s,\bm{\lambda})=\partial^2 H_{t_0}/\partial s^2 (s,\bm{\lambda})=\partial^3 H_{t_0}/\partial s^3 (s,\bm{\lambda})=\partial^4 H_{t_0}/\partial s^4 (s,\bm{\lambda})=0$ if and only if $\kappa _g(s,t_0)\pm \kappa _n(s,t_0)\not= 0$, $((\kappa _n\pm\kappa_g)\tau _g\mp (\kappa _n'\pm\kappa _g'))(s_0,t_0)=((\kappa _n\pm\kappa_g)\tau _g\mp (\kappa _n'\pm\kappa _g'))'(s,t_0)=0$  and 
 \[
 \bm{\lambda}=\bm{\Gamma} (s,t_0)+\displaystyle{\frac{1}{\kappa _g(s,t_0)\pm \kappa _n(s,t_0)}}(\bm{b}(s,t_0)\pm \bm{n}(s,t_0)).\]
\end{Pro}
\demo
Since we have the pseudo-orthonormal frame $\{\bm{\Gamma}(s,t),\bb (s,t), \bn (s,t), \bt(s,t)\},$
there exist real numbers $\lambda, \mu, \nu \in \R$ such that $\bm{\lambda}=\xi\bm{\Gamma}(s,t)+
\lambda \bb (s,t_0)+\mu \bn(s,t_0) +\nu \bt(s,t_0).$
\par\noindent
(1) The condition  $\partial H_{t_0}/\partial s (s,\bm{\lambda})=0$ means that $\nu =0$.
Moreover, the condition $H_{t_0}(s,\bm{x})=0$ means that $\xi=1$. 
Since $\langle \bm{\lambda},\bm{\lambda}\rangle =-1,$
 we have $\lambda ^2-\mu ^2=0.$ It follows that
\[
\bm{\lambda}=\bm{\Gamma}(s,t_0)+
\mu (\bb (s,t_0)\pm \bn(s,t_0)).
\]
We put $u=\mu $. This completes the proof of (1).
\par\noindent 
(2) With the assumption that (1) holds, the condition $\partial^2 H_{t_0}/\partial s^2 (s,\bm{\lambda})=0$
means that 
\[
0=\langle \bm{\Gamma}-\kappa _g\bm{b}+\kappa _n\bm{n},\bm{\lambda}\rangle=(\kappa _g\pm\kappa _n)u-1.
\]
Therefore, we have $\kappa _g(s,t_0)\pm \kappa _n(s,t_0)\not= 0$ and \[
 \bm{\lambda}=\bm{\Gamma} (s,t_0)+\displaystyle{\frac{1}{\kappa _g(s,t_0)\pm \kappa _n(s,t_0)}}(\bm{b}(s,t_0)\pm \bm{n}(s,t_0)).
 \]
This completes the proof of (2).
\par\noindent
(3) By the similar arguments to the above cases, we have the assertion (3).
\par
Moreover, if we calculate the $4$th derivative $\displaystyle{\frac{\partial^4 H_{t_0}}{\partial s^4}}$, then
we have the assertion (4). Since those arguments are tedious, we omit the detail here.
\enD

According to the above proposition, we introduce an invariant defined by
\[
\sigma^\pm (s,t)=((\kappa _n\pm\kappa _g)\tau _g\mp(\kappa_n'\pm\kappa_g'))(s,t).
\]
\begin{Pro} Suppose that $\kappa _g(s,t_0)\pm \kappa _n(s,t_0)\not= 0$ and we denote $\tau =+\ \mbox{or}\ -.$  Then the following conditions are equivalent{\rm :}
\par
{\rm (1)} $\sigma ^\tau (s,t_0)\equiv 0$,
\par
{\rm (2)}  $\{\bm{\lambda}^\tau_0\}=\mathbb{LF}^\tau_{\mathcal{S}_{t_0}}$,
\par
{\rm (3)} There exists $\bm{\lambda}_0\in AdS^3$ such that $\mathcal{S}_{t_0}\subset LC^{AdS}(\bm{\lambda}_0).$
\end{Pro}
\demo
We define $\bm{\ell }_\pm:I\lon AdS^3$ by
\[
\bm{\ell}_\pm (s)=\bm{\Gamma} (s,t_0)+\displaystyle{\frac{1}{\kappa _g(s,t_0)\pm \kappa _n(s,t_0)}}(\bm{b}(s,t_0)\pm \bm{n}(s,t_0)).
\]
Then $\bm{\ell}_\pm (I)=\mathbb{LF}^\pm_{\mathcal{S}_{t_0}}.$
By a straightforward calculation, we have 
\[
\bm{\ell}_\pm '(s)=-\frac{\sigma ^\pm (s,t_0)}{(\kappa _g(s,t_0)\pm \kappa _n(s,t_0))^2}(\bn(s,t_0)\pm \bb(s,t_0)).
\]
Therefore conditions (1) and (2) are equivalent. Suppose that (2) holds.
Then we have $\bm{\lambda}_0^\tau=\bm{\ell}_\tau (s)$ for any $s\in I.$
Thus, we have $\bm{\Gamma} (s,t_0)\in \Lambda _{\bm{\lambda}^\tau_0}\cap AdS^3=LC^{AdS}({\bm{\lambda}^\tau_0})$ for any $s\in I$, so that  (3) holds.
Suppose that (3) holds. Then there exists a point $\bm{\lambda}_0\in AdS^3$ such that
$\mathcal{S}_{t_0}\subset LC^{AdS}({\bm{\lambda}_0})=HP(\bm{\lambda}_0,-1)\cap AdS^3.$
This condition is equivalent to the condition that
$\langle \bm{\Gamma}(s,t_0),\bm{\lambda}_0\rangle =-1$  at any $s\in I.$
Then $H_{t_0}(s,\bm{\lambda}_0)$ is constantly equal to zero.
By the previous calculations, this is equivalent to the condition that $\{\bm{\lambda} _0 \}=\bm{\ell}_\tau (I)$ and (1) holds.
This completes the proof.
\enD
\par
We also have a classification of singularities of momentary lightlike hypersrufaces.
\begin{Th} {\rm (1)} The lightlike hypersruface $\mathbb{LS}^\pm_{\mathcal{S}_{t_0}}(I\times\{t_0\}\times \R)$ at $\bm{\lambda}_0=\bm{\ell}_\pm (s_0)\in \mathbb{LF}^\pm_{\mathcal{S}_{t_0}}$ is local diffeomorphic to the cuspidaledge $\bm{CE}$ if $\sigma ^\pm (s_0,t_0)\not= 0,$
\par
{\rm (1)} The lightlike hypersruface $\mathbb{LS}^\pm_{\mathcal{S}_{t_0}}(I\times\{t_0\}\times \R)$ at $\bm{\lambda}_0=\bm{\ell}_\pm (s_0)\in \mathbb{LF}^\pm_{\mathcal{S}_{t_0}}$ is local diffeomorphic to the swallowtail $\bm{SW}$ if $\sigma ^\pm (s_0,t_0)= 0$ and $\partial \sigma ^\pm/\partial s (s_0,t_0)\not= 0.$
\par
Here, $\bm{CE}=\{(u,v^2,v^3)\in (\R^3,0)\ |\ (u,v)\in (\R^2,0)\ \}$ and $\bm{SW}=\{(3u^4+vu^2,4u^2+2uv,v)\in (\R^3,0)\ |\  (u,v)\in (\R^2,0)\ \}.$
\end{Th}

\par
In order to prove Theorem 10.3, we use some general results on the singularity theory for unfoldings of function germs.
Detailed descriptions are found in the book \cite{Bru-Gib}.
Let $F:({\R}\times{\R}^r,(s_0,x_0))\rightarrow {\R}$ be a function germ.
We call $F$ an {\it $r$-parameter unfolding\/} of $f$, where $f(s)=F_{x_0}(s,x_0).$
We say that $f$ has an {\it $A_k$-singularity\/} at $s_0$ if $f^{(p)}(s_0)=0$ for all $1\leq p\leq k$,
and $f^{(k+1)}(s_0)\ne 0.$
Let $F$ be an unfolding of $f$ and $f(s)$ has an $A_k$-singularity $(k\geq 1)$ at $s_0.$
We denote the $(k-1)$-jet of the partial derivative
$\frac{\partial F}{\partial x_i}$ at $s_0$ by 
$j^{(k-1)}(\frac{\partial F}{\partial x_i}(s,x_0))(s_0)=\sum_{j=0}^{k-1} \alpha_{ji}(s-s_0)^j$ 
for $i=1,\dots ,r$. Then
$F$ is called an {\it $\mathcal{R}$-versal unfolding \/} if the $k\times{r}$ 
matrix of coefficients $(\alpha _{ji})_{j=0,\dots ,k-1;i=1,\dots ,r}$
has rank $k$ $(k\leq {r}).$
We introduce an important set concerning the unfoldings relative to the above notions. 
A {\it $\ell$th-discriminant set} of $F$ is 
$$
{\mathcal D}^\ell _F=\left\{x\in {\R}^r\Bigm|\exists s \ {\rm with }\ F=
\frac{\partial F}{\partial s}=\cdots =\frac{\partial^\ell F}{\partial s^\ell}= 0 \ {\rm at }\ (s,x)\right\}.$$
For $\ell =1,$ it is simply denoted by $\mathcal{D}_F,$ which is called a {\it discriminant set} of $F.$
Then we have the following classification (cf., \cite{Bru-Gib}).
\begin{Th} Let $F:({\R}\times{\R}^r,(s_0,x_0))\rightarrow {\R}$
be an $r$-parameter unfolding of $f(s)$ which has an $A_k$ singularity at $s_0$.
Suppose that $F$ is an $\mathcal{R}$-versal unfolding.
\par
{\rm (1)} If $k=2$, then ${\mathcal D}_F$ is locally diffeomorphic to  
$\bm{CE}\times {\R}^{r-2}$.
\par
{\rm (2)} If $k=3$, then ${\mathcal D}_F$ is locally diffeomorphic to 
$\bm{SW}\times {\R}^{r-2}$.
\end{Th}
\par
For the proof of Proposition 10.3, we have the following propositions.
Let  $\bm{\Gamma} :I\times J\lon W\subset \R^3_1$ be a world sheet with  $\kappa _n(s,t)\pm \kappa _g(s,t)\not= 0$
and $H:I\times J\times \R^3\lon \R$ the $AdS$-height function on $\bm{\Gamma}.$
We define $h_{t_0,\bm{\lambda}_0}(s)=H_{t_0}(s,\bm{\lambda}_0)=H(s,t_0,\bm{\lambda}_0)$ and consider that
$H_{t_0}$ is a $3$-parameter unfolding of $h_{t_0,\bm{\lambda}_0}.$
\begin{Pro}   If $h_{t_0,\bm{\lambda}_0}$ has an $A_k$-singularity $(k=2,3)$ at $s_0,$
then $H_{t_0}$ is an $\mathcal{R}$-versal unfolding of $h_{t_0,\bm{\lambda}_0}.$
\end{Pro}
\demo
We write that $\bm{\Gamma}(s,t)=(X_0(s,t),X_1(s,t),X_2(s,t))$ and $\bm{\lambda}=(\lambda_{-1},\lambda _0,\lambda _1,\lambda _2).$
Then we have
\[
H_{t_0}(s,\bm{\lambda}_0)=-X_{-1}(s,t_0)\lambda_{-1} -X_0(s,t_0)\lambda _0+X_1(s,t_0)\lambda _1+X_2(s,t_0)\lambda _2+1.
\]
Since $\bm{\lambda}\in AdS^3,$ we have
$
-\lambda ^2_{-1}-\lambda ^2_0+\lambda ^2_1+\lambda ^2_2=-1.
$
Then we consider the local coordinates $(\lambda _0,\lambda _1,\lambda_2)$ of $AdS^3$ given by $\lambda _{-1}=\sqrt{1-\lambda^2_0+\lambda ^2_1+\lambda ^2_2}>0.$
Therefore, we have
\[
\frac{\partial H_{t_0}}{\partial \lambda _{0}}(s,\bm{\lambda}_0)=-X_0(s,t_0)+X_{-1}(s,t_0)\frac{\lambda_0}{\lambda _{-1}},\  \frac{\partial H_{t_0}}{\partial \lambda_i}(s,\bm{\lambda}_0)=X_i(s,t_0)-X_{-1}(s,t_0)\frac{\lambda_i}{\lambda _{-1}},\ i=1,2.
\]
Thus we obtain
\begin{eqnarray*}
j^2\left(\frac{\partial H_{t_0}}{\partial \lambda_0}(s_0,\bm{\lambda}_0)\right)&{}&\!\!\!\!\!\!=-X_0(s_0,t_0)+X_{-1}(s_0,t_0)\frac{\lambda_0}{\lambda _{-1}}
\\
&{}&+\left(-\frac{\partial X_0}{\partial s}(s_0,t_0)+\frac{\partial X_{-1}}{\partial s}(s_0,t_0)\frac{\lambda_0}{\lambda _{-1}}\right)(s-s_0)\\
&{}&+\frac{1}{2}\left(-\frac{\partial^2 X_0}{\partial s^2}(s_0,t_0)+\frac{\partial^2 X_{-1}}{\partial s^2}(s_0,t_0)\frac{\lambda_0}{\lambda _{-1}}\right)(s-s_0)^2,
\end{eqnarray*}
\begin{eqnarray*}
j^2\left(\frac{\partial H_{t_0}}{\partial \lambda_i}(s_0,\bm{\lambda}_0)\right)&{}&\!\!\!\!\!\!=X_i(s_0,t_0)-X_{-1}(s_0,t_0)\frac{\lambda_i}{\lambda _{-1}}\\
&{}&+\left(\frac{\partial X_i}{\partial s}(s_0,t_0)-\frac{\partial X_{-1}}{\partial s}(s_0,t_0)\frac{\lambda_i}{\lambda _{-1}}\right)(s-s_0)\\
&{}&+\frac{1}{2}\left(\frac{\partial^2 X_i}{\partial s^2}-\frac{\partial^2 X_{-1}}{\partial s^2}(s_0,t_0)\frac{\lambda_i}{\lambda _{-1}}\right)(s-s_0)^2,
\end{eqnarray*}
$i=1,2.$
We consider a matrix
\[
\mbox{\Large $A$}=\begin{pmatrix} -X_0+X_{-1}\frac{\lambda_0}{\lambda _{-1}} & X_1-X_{-1}\frac{\lambda_0}{\lambda _{-1}} & X_2-X_{-1}\frac{\lambda_0}{\lambda _{-1}}\\
 -\frac{\partial X_0}{\partial s}+\frac{\partial X_{-1}}{\partial s}\frac{\lambda_0}{\lambda _{-1}} & \frac{\partial X_1}{\partial s}-\frac{\partial X_{-1}}{\partial s}\frac{\lambda_0}{\lambda _{-1}} &
\frac{\partial X_2}{\partial s}-\frac{\partial X_{-1}}{\partial s}\frac{\lambda_0}{\lambda _{-1}} \\
 -\frac{\partial^2 X_0}{\partial s^2}+\frac{\partial^2 X_{-1}}{\partial s^2}\frac{\lambda_0}{\lambda _{-1}} & \frac{\partial^2 X_1}{\partial s^2}-\frac{\partial^2 X_{-1}}{\partial s^2}\frac{\lambda_1}{\lambda _{-1}} &
\frac{\partial^2 X_2}{\partial s^2}-\frac{\partial^2 X_{-1}}{\partial s^2}\frac{\lambda_2}{\lambda _{-1}}
\end{pmatrix}
\]
at $(s_0,t_0).$
Then we have
\[
\det \mbox{\Large $A$}=\frac{1}{\lambda_{-1}}\left\langle \bm{\lambda}_0,\bm{\Gamma}(s_0,t_0)\wedge \frac{\partial \bm{\Gamma}}{\partial s}(s_0,t_0)\wedge 
\frac{\partial^2 \bm{\Gamma}}{\partial s^2}(s_0,t_0)\right\rangle
\]
We also have 
\[
\frac{\partial \bm{\Gamma}}{\partial s}(s_0,t_0)=\bm{t}(s_0,t_0),\ 
\frac{\partial^2 \bm{\Gamma}}{\partial s^2}(s_0,t_0)=-\kappa _g(s_0,t_0)\bm{b}(s_0,t_0)+\kappa_n(s_0,t_0)\bm{n}(s_0,t_0).
\]
By Proposition 10.1, we have $\bm{\lambda}_0=(\bm{\Gamma}+(\bm{b}\pm\bm{n})/(\kappa_g\pm\kappa_n))(s_0,t_0)$, so that 
\[
\det \mbox{\Large $A$}=\frac{1}{\lambda_{-1}}\langle \bm{\lambda}_0,\kappa _g(s_0,t_0)\bm{n}(s_0,t_0)-\kappa _n\bm{b}(s_0,t_0)\rangle
=\pm\frac{1}{\lambda_{-1}}\not=0.
\]
 This means that $H_{t_0}$ is an $\mathcal{R}$-versal unfolding of $h_{t_0,\bm{\lambda}_0}.$
 \par
 For other local coordinates of $AdS^3,$ we have the similar calculations to the above case.
 \enD
\vskip3pt
\par\noindent
{\it Proof of Theorem 10.3.}
By (1) of Proposition 10.1, the discriminant set $D_{H_{t_0}}$ of the $AdS$-height
function on $\mathcal{S}_{t_0}$ is the lightlike hypersruface along $\mathcal{S}_{t_0}.$
It also follows  (3) and (4) of Proposition 10.1 that
$h_{t_0,\bm{\lambda}_0}$ has an $A_2$-singularity (respectively, $A_3$-singularity) at $s_0$ if $\sigma ^\pm (s_0,t_0)\not= 0$ (respectively, $\sigma ^\pm (s_0,t_0)= 0$ and $(\sigma ^\pm)' (s_0,t_0)\not= 0$).
By Proposition 10.5, $H_{t_0}$ is an $\mathcal{R}$-versal unfolding of $h_{t_0,\bm{\lambda}_0}$ for each case.
Then we can apply the classification theorem (Theorem 10.4) to our situation.
This completes the proof.
\enD
\par We remark that $D^2_{H_{t_0}}$ is the lightlike focal curve $\mathbb{LF}^\pm_{\mathcal{S}_{t_0}}$.
Since the critical value set of the swallow tail is locally diffeomorphic to a {\it $(2,3,4)$-cusp}
which is defined by
$C=\{(t^2,t^3,t^4)\ |\ t\in \R\},$ we have the following corollary.
\begin{Co}
The lightlike focal curve $\mathbb{LF}^\pm_{\mathcal{S}_{t_0}}$ is locally diffeomorphic to a line if $\sigma ^\pm (s_0,t_0)\not= 0$.
It is locally diffeomorphic to the $(2,3,4)$-cusp if $\sigma ^\pm (s_0,t_0)= 0$ and $(\sigma ^\pm)' (s_0,t_0)\not= 0$.
\end{Co}  
\par
On the other hand, we now classify $S.P^+$-Legendrian stable graph-like Legendrian unfoldings $\mathscr{L}_H(\Sigma _*(H))$ by
$S.P^+$-Legendrian equivalence. By Theorems 5.5 and 5.6, it is enough to classify $\overline{f}$ by $S.P$-$\mathcal{K}$-equivalence under the condition that
\[
\dim _{\R} \frac{{\cal E}_{1+1}}{\left\langle \frac{\partial \overline{f}}{\partial q},\overline{f} \right\rangle _{{\cal E}_{1+1}}+
\left\langle \frac{\partial \overline{f}}{\partial t} \right\rangle _{\R}} \leq 3.
\]
In \cite{Izudoc,Izu95} we have the following proposition.
\begin{Pro}
With the above condition, $\overline{f}:(\R\times\R,0)\lon (\R,0)$ with $\partial\overline{f}/\partial t(0)\not= 0$ is $S.P$-$\mathcal{K}$-equivalent to one of the following germs{\rm :}
\par
{\rm (1)} $q,$
\par
{\rm (2)} $\pm t\pm q^2,$
\par
{\rm (3)} $\pm t +  q^3,$
\par
{\rm (4)} $\pm t \pm q^4,$
\par
{\rm (5)} $\pm t + q^5.$
\end{Pro}

The infinitesimally $S.P^+$-$\mathcal{K}$-versal unfolding
$\mathcal{F}:(\R\times (\R^3\times \R),0)\lon (\R,0)$ of each germ in the above list is given as follows (cf. \cite[Theorem 4.2]{Izu95}):
\par
(1) $q$
\par
(2) $\pm t\pm q^2,$
\par
(3) $\pm t + q^3+x_0q,$
\par
(4) $\pm t\pm q^4+x_0q+x_1q^2,$
\par
(5) $\pm t+q^5+x_0q+x_1q^2+x_2q^3.$
\par\noindent

By Theorem 5.6, we have the following classification.
\begin{Th} Let $(W,\mathcal{S})$ be a world sheet in $AdS^3$ parametrized by a timelike embedding $\bm{\Gamma} :J\times I\lon AdS^3$
and $H:J\times I\times AdS^3\lon \R$ be the $AdS$-height squared function of $(W,\mathcal{S}).$
Suppose that the corresponding graph-like Legendrian unfolding $\mathscr{L}_H(\Sigma _*(H))\subset J^1(AdS^3,I)$ is $S.P^+$-Legendrian stable.
Then the germ of the image of the unfolded lightlike hypersrufaces $\mathbb{LH}_W$ at any point is $S.P^+$-diffeomorphic to one of the following set germs in 
$(\R^3\times \R,0)${\rm :}
\par
{\rm (1)} $\{(u,v,w),0)\ |\ (u,v,w)\in (\R^3,0)\ \},$

\par
{\rm (2)} $\{(-u^2,v,w),\pm 2u^3)\ |\ (u,v,w)\in (\R^3,0)\ \},$
\par
{\rm (3)} $\{(\mp 4u^3-2vu,v,w),3u^3\pm vu^2)\ |\ (u,v,w)\in (\R^3,0)\ \},$
\par
{\rm (4)} $\{((5u^4+2vu+3wu^2,v,w),\pm(4u^4+vu^2+2wu^3))\ |\ (u,v,w)\in (\R^3,0)\ \}.$
\end{Th}
\demo
For any $(s_0,t_0,\bm{\lambda}_0)\in J\times I\times AdS^3,$ the germ of $\mathscr{L}_H(\Sigma _*(GH)\subset J^1(AdS^3,I)$ at $\bm{z}_0=\mathscr{L}_H(s_0,t_0,\bm{\lambda}_0)$
is $S.P^+$-Legendrian stable. It follows that the germ of $h_{\bm{\lambda}_0}$ at $(s_0,t_0)$ is $S.P$-$\mathcal{K}$-equivalent to one of the germs in the list of Proposition 10.7.
By Theorem 5.6, the graph-like Legendrian unfolding $\mathscr{L}_H(\Sigma _*(H))$ is $S.P^+$-Legendrian equivalent to the graph-like Legendrian unfolding 
$\mathscr{L}_{\mathcal{F}}(\Sigma_*(\mathcal{F}))$ where $\mathcal{F}$ is the infinitesimally $S.P$-$\mathcal{K}$-versal unfolding of one of 
the germs in the list of Proposition 10.7.
It is also equivalent to the condition that the germ of the graph-like big front $W(\mathscr{L}_{\mathcal{F}}(\Sigma_*(\mathcal{F})))$ is $S.P^+$-diffeomorphic to the corresponding graph-like big front of one of the normal forms.
For each normal form, we can obtain the graph-like big front.
We only show that (5) in Proposition 10.7.
In this case we consider $\mathcal{F}(q,x_0,x_1,x_2,t)=\pm t+q^5+x_0q+x_1q^2+x_2q^3.$
Then we have
\[
\frac{\partial \mathcal{F}}{\partial q}=5q^4+x_0+2x_1q+3x_1q^2,
\]
so that the condition $\mathcal{F}=\partial \mathcal{F}/\partial q=0$ is equivalent to the condition that
\[
x_0=-(5q^4+x_0+2x_1q+3x_1q^2),\ t_0=\pm(4q^5+x_1q^2+2x_2q^3).
\]
If we put $u=q, v=x_0,w=x_1,$ then we have
\[
W(\mathscr{L}_{\mathcal{F}}(\Sigma_*(\mathcal{F})))=\{((-(5u^4+2vu+3wu^2),v,w),\pm(4u^4+vu^2+2wu^3))| (u,v,w)\in (\R^3,0)\}.
\]
It is $S.P^+$-diffeomorphic to the set germ of (4).
We have similar calculations for other cases.
We only remark here that we obtain the germ of (1) for both the germs of (1) and (2) in Proposition 10.7.
Since $W(\mathscr{L}_{\mathcal{H}}(\Sigma_*(\mathcal{H})))=\mathbb{LH}_W,$ this completes the proof.
\enD
As a corollary, we have a local classification of BR-caustics in this case.
\begin{Co} With the same assumption for the world sheet $(W,\mathcal{S})$ as Theorem 10.8, the BR-caustic $C(W,\mathcal{S})$ of $(W,\mathcal{S})$ at a singular point is locally diffeomorphic to the cuspidaledge $\bm{CE}$ or the swallowtail $\bm{SW}$.
\end{Co}
\demo
The BR-caustic $C(W,\mathcal{S})$ of $(W,\mathcal{S})$ is the set of the critical values of $\pi_1\circ \overline{\pi}|_{\mathscr{L}_H(\Sigma _*(H))}.$
Therefore, it is enough to calculate the set of critical values of $\pi_1\circ \overline{\pi}|_{\mathscr{L}_{\mathcal{F}}(\Sigma_*(\mathcal{F}))}$
for each normal form $\mathcal{F}$ in Proposition 10.7.
For (5) in Proposition 10.7, by the proof of Theorem 10.8 we have
$$
\Sigma_*(\mathcal{F})=\{(u,5u^4+2vu+3wu^2,v,w)\in (\R\times (\R^3\times \R),0)|(u,v,w)\in (\R^3,0)\}.
$$
It follows that
$$
\pi_1\circ \overline{\pi}\circ \mathscr{L}_{\mathcal{F}}(u,5u^4+2vu+3wu^2,v,w)=(5u^4+2vu+3wu^2,v,w).
$$
Then the Jacobi matrix of $f(u,v,w)=(5u^4+2vu+3wu^2,v,w)$ is
\[
J_f=\begin{pmatrix}
20u^3+2v+6wu & 0 & 0 \cr
2u & 1 & 0 \cr
3u^2 & 0 & 1 \cr
\end{pmatrix},
\]
so that the set of critical values of $f$ is  given by
\[
\{(-(15u^4+3wu^2),-10u^3-3wu,w)\in (\R^3,0)\ |\ (u,w)\in (\R^2,0)\}.
\] 
For a linear isomorphism $\psi:(\R^3,0)\lon \R^3,0)$ defined by
$
\psi (x_0,x_1,x_2)=(-\frac{1}{5}x_0,-\frac{2}{5}x_1,\frac{3}{5}x_2),
$
we have
$\psi(-(15u^4+3wu^2),-10u^3-3uw,w)=(3u^4+\frac{3}{5}wu^2,4u^3+\frac{6}{5}wu,\frac{3}{5}w).$
If we put $U=u, V=\frac{3}{5}w,$ then we have $(3U^4+VU^2,4U^3+2VU,V),$ which is the parametrization of $\bm{SW}.$
By the arguments similar to the above, we can show that the set of critical values of $\pi_1\circ \overline{\pi}|_{\mathscr{L}_{\mathcal{F}}(\Sigma_*(\mathcal{F}))}$ is a regular surface for (3) and is diffeomorphic to $\bm{CE}$ for (4) in Proposition 10.7, respectively.
This completes the proof.
\enD
\begin{Rem}{\rm
Since a world sheet $(W,\mathcal{S})$ is a timelike surface in $AdS^3,$ we can define the {\it $AdS$-evolute} of $(W,\mathcal{S})$ by
\[
Ev^{AdS}_{(W,\mathcal{S})}=\displaystyle{\bigcup _{i=1}^2\left\{\frac{\pm 1}{\sqrt{\kappa ^2_i(u,t)-1}}(\kappa _i(u,t)\bm{X}(u,t)+\bm{n}^S(u,t))\ |\ (u,t)\in U\times I,\kappa ^2_i(u,t)>1\ \right\}},
\]
where $\kappa _i(s,t)$ $(i=1,2)$ are the principal curvatures of $W$ at $p=\bm{X}(u,t)$ with respect to $\bm{n}^S$ (cf. \cite{timeads}).
The $AdS$-evolute of a timelike surface has singularities in general. 
Actually, it is a caustic in the the theory of Lagrangian singularities.
Similar to the notion of evolutes of surfaces in Euclidean space $\R^3$ (cf. \cite{Porteous}), the corank two singularities of the $AdS$-evolute appear at the umbilical points (i.e. $\kappa _1(u,t)=\kappa _2(u,t)$).
The singularities of the $AdS$-evolute of a generic surface in $AdS^3$ are classified into $\bm{CE}$, $\bm{SW},$ $\bm{PY}$ or $\bm{PU}$,
where $\bm{PY}=\{(u^2-v^2+2uv,-2uv+2uw,w)|w^2=u^2+v^2\}$ is the {\it pyramid} and $\bm{PU}=\{(3u^2+wv,3v^2+wu,w)|w^2=36uv\}$ is the {\it purse}.
The pyramid and the purse of the $AdS$-evolute correspond to the umbilical points of the timelike surface in $AdS^3.$
So the singularities of BR-caustics of world sheets are different from those of the $AdS$-evolutes of surfaces.
Since the singularities of BR-caustics are only corank one singularities, the pyramid and the purse never appeared in general.
Moreover, the normal geodesic of a timelike surface is a spacelike curve, so that it is not a ray in the sense of the relativity theory.
Therefore, the $AdS$-evolute of a timelike surface in anti-de Sitter space-time is not a caustic in the sense of physics. 
}
\end{Rem}

{\small
\par\noindent
Shyuichi Izumiya, Department of Mathematics, Hokkaido University, Sapporo 060-0810,Japan
\par\noindent
e-mail:{\tt izumiya@math.sci.hokudai.ac.jp}
}

\begin{thebibliography}{99999}
{\renewcommand{\baselinestretch}{1}
\small
\bibitem{Arnold1}
V. I. Arnol'd, S. M. Gusein-Zade and A. N. Varchenko,
\newblock{\em Singularities of Differentiable Maps vol. I}.
\newblock Birkh\"auser, 1986.

\bibitem{Arnold-pink}
V. I. Arnol'd,
\newblock{\em Singularities of caustics and wave fronts}.
\newblock Math. Appl. 62, Kluwer , Dordrecht, 1990.


\bibitem{Bousso}
R. Bousso,
\newblock{\em The holographic principle},
\newblock REVIEWS OF MODERN PHYSICS. 74 (2002), 825--874.

\bibitem{Bousso-Randall}
R. Bousso  and L Randall, 
\newblock{\em Holographic domains of ant-de Sitter space}, 
\newblock Journal of High Energy Physics. 04 (2002), 057

\bibitem{Brekhov}
L. Brekhovskikh,
\newblock{\em Wave in layered media},
\newblock Academic press, 1980


\bibitem{Bru-Gib}
J. W. Bruce and P. J. Giblin,
\newblock{\em Curves and singularities (second edition)}.
\newblock Cambridge University press (1992)


\bibitem{Ch}
S. Chandrasekhar,
\newblock{\em The Mathematical Theory of Black Holes},
International Series of Monographs on Physics. 69
\newblock Oxford University press, 1983.


\bibitem{timeads}
L. Chen, S. Izumiya and D. Pei,
\newblock{\em Timelike hypersurfaces in anti-de Sitter space from a contact viewpoint},
Journal of Mathematical Sciences, {\bf 199} (2014), 629--645


\bibitem{Hormander}
L. H\"ormander,
\newblock{\em Fourier Integral Operators,I}.
\newblock Acta. Math. {\bf 128} (1972), 79--183


\bibitem{Izudoc}
S. Izumiya, 
\newblock{\em Generic bifurcations of varieties}.
\newblock manuscripta math. {\bf 46} (1984), 137--164


\bibitem{Izumiya93}
S. Izumiya,
\newblock{\em Perestroikas of optical wave fronts and graphlike Legendrian unfoldings}.
\newblock J. Differential Geom. {\bf 38} (1993), 485--500.

\bibitem{Izu95}
S. Izumiya,
\newblock{\em Completely integrable holonomic systems of first-order differential equations}.
\newblock Proc. Royal Soc. Edinburgh 125A (1995), 567--586.

\bibitem{IzuHJ93}
S. Izumiya,
\newblock{\em Geometric singularities for Hamilton-Jacobi equations},
\newblock Adv. Studies in Pure Math. {\bf 22} (1993), 89--100

\bibitem{Izu07}
S. Izumiya,
\newblock{\em Differential Geometry from the viewpoint of Lagrangian or Legendrian singularity theory}.
\newblock in {\it Singularity Theory, Proceedings of the 2005 Marseille Singularity School and Conference, by D. Ch\'eniot et al}. World Scientific (2007) 241--275.


\bibitem{IzuSM} S. Izumiya and M.C. Romero Fuster, {\em The lightlike flat geometry
on spacelike submanifolds of codimension two in Minkowski space.} Selecta Math. (N.S.) 13 (2007), no. 1, 23--55.


\bibitem{Izumiya-Takahashi} S. Izumiya and M. Takahashi,
\newblock{\em Spacelike parallels and evolutes in Minkowski pseudo-spheres}.
\newblock{Journal of Geometry and Physics}. {\bf 57} (2007), 1569--1600.

\bibitem{Izumiya-Takahashi2} S. Izumiya and M. Takahashi,
\newblock{\em Caustics and wave front propagations: Applications to differential geometry}.
\newblock{Banach Center Publications. Geometry and topology of caustics}. {\bf 82} (2008)  125--142.

\bibitem{Izumiya-Takahashi3} S. Izumiya and M. Takahashi,
\newblock{\em Pedal foliations and Gauss maps of hypersurfaces in Euclidean space}.
\newblock{Journal of Singularities}. {\bf 6} (2012)  84--97.

\bibitem{Graph-like}
S. Izumiya, \newblock{\em The theory of graph-like Legendrian unfoldings and its applications},
 J. of Singularities, {\bf 12} (2015), 53--79 DOI:10.5427/jsing.2015.12d
 
\bibitem{IzuLight} 
S. Izumiya, \newblock{\em Lightlike hypersurfaces along spacelike submanifolds in anti-de Sitter space},
preprint (2013)

\bibitem{GeomLag14}
S. Izumiya, \newblock{\em Geometric interpretation of Lagrangian equivalence},
preprint (2014).

\bibitem{KR01}
A. Karch and L. Randall,
{\em Locally localized gravity.}
J. High Energy Physics {\bf 05}, (2001), 008.

\bibitem{M98}
M. Maldacena,
{\rm The Large N Limit of Superconformal Field Theories and
Supergravity}, Adv. Theor. Math. Phys. {\bf 2} (1998) 231--252.


\bibitem{mont1}
J. A. Montaldi,
\newblock {\em On contact between submanifolds},
\newblock Michigan Math. J., {\bf 33} (1986), 81--85.

\bibitem{MTW}
C. W. Misner, K. S. Thorne and J. W. Wheeler,
\newblock{\em Gravitation},
\newblock W. H. Freeman and Co., San Francisco, CA (1973).


\bibitem{Oneil}
B. O'Neill,
\newblock{\em Semi-Riemannian Geometry},
\newblock Academic Press, New York, 1983.


\bibitem{Nye}
J. F. Nye,
\newblock{\em Natural focusing and fine structure of light}.
\newblock Institute of Physics Publishing, Bristol and Philadelphia, 1999


\bibitem{Peters}
A. O. Petters, H. Levine and J. Wambsganss,
\newblock{\em Singularity theory and gravitational lensing},
\newblock Birkh\"auser, 2001



\bibitem{Penrose}
R. Penrose,
\newblock{\em Null Hypersurface Initial Data for Classical Fields of Arbitrary Spin and for General Relativity},
General Relativity and Gravitation, {\bf 12} (1963), 225--264

\bibitem{Porteous}
I. Porteous,
\newblock{\em The normal singularities of submanifold}.
\newblock J. Diff. Geom. {5}, (1971), 543--564.


 \bibitem{RS99}
 L. Randall and R. Sundrum,
 {\em An alternative to Compactification},
 Physical Review Letters, {\bf 83} (1999), 4690--4693.

\bibitem{Saloom}
A. Saloom and F. Tari,
\newblock{\em Curves in the Minkowski plane and their contact with pseudo-circles},
Geom. Dedicata, {\bf 159}, (2012), 109--124 

\bibitem{Siino}
M. Siino and T. Koike,
\newblock{\em Topological classification of black holes: generic Maxwell set and crease set of a horizon},
International Journal of Modern Physics D, {\bf 20}, (2011), 1095.

\bibitem{Farid}
F. Tari,
\newblock{\em Caustics of surfaces in the Minkowski $3$-space},
Quarterly Journal of Mathematics, {\bf 63}, (2012), 189--209.


\bibitem{W98}
E. Witten, {\em Anti de Sitter space and holography},
Adv. Theor. Math. Phys. {\bf 2} (1998), 253--291.

\bibitem{Zak}
V. M. Zakalyukin,
\newblock{\em Lagrangian and Legendrian singularities},
\newblock Funct. Anal. Appl. (1976), 23--31.


\bibitem{Zak1}
V. M. Zakalyukin,
\newblock{\em Reconstructions of fronts and caustics depending one
parameter and versality of mappings},
\newblock J. Sov. Math. 27 (1984), 2713--2735.


\bibitem{Zakalyukin95} 
V.M. Zakalyukin, 
\newblock{\em Envelope of Families of Wave Fronts and Control Theory.}
\newblock Proc. Steklov Inst. Math. {\bf 209} (1995), 114--123.

}
\end{thebibliography}
 \end{document}